\documentclass[oneside, 12pt]{amsart}
\usepackage{amscd, amssymb, amsmath, mathrsfs}
\usepackage[english]{babel}
\usepackage{url}
\usepackage{tikz-cd}
\usepackage{booktabs}
\usepackage[pdftex, colorlinks=true,  citecolor=blue, linkcolor=blue, linktocpage=true]{hyperref}

\DeclareUrlCommand\arXiv{\urlstyle{same}}

\newcommand\mylabel[1]{\label{#1}\marginpar{\vspace{-1ex}\medskip\medskip\footnotesize \tt #1}}
\renewcommand\mylabel[1]{\label{#1}}
\newcommand{\mydate}{
\number\day\space
\ifcase\month \or January\or February\or March\or April\or May\or June\or July\or August\or September\or October\or November\or December\fi 
\space\number\year}

\newtheorem{theorem}{Theorem}[section]
\newtheorem*{maintheorem}{Theorem}
\newtheorem{lemma}[theorem]{Lemma}
\newtheorem{proposition}[theorem]{Proposition}
\newtheorem{corollary}[theorem]{Corollary}

\theoremstyle{definition}
\newtheorem{definition}[theorem]{Definition}

\newtheorem*{acknowledgement}{Acknowledgement}

\theoremstyle{remark}

\newcommand{\ZZ}{\mathbb{Z}}

\newcommand{\FF}{\mathbb{F}}

\newcommand{\PP}{\mathbb{P}}
\renewcommand{\AA}{\mathbb{A}}
\newcommand{\GG}{\mathbb{G}}

\newcommand{\ideala}{\mathfrak{a}}

\newcommand{\idealh}{\mathfrak{h}}

\newcommand{\shA}{\mathscr{A}}

\newcommand{\shE}{\mathscr{E}}
\newcommand{\shF}{\mathscr{F}}

\newcommand{\shN}{\mathscr{N}}
\newcommand{\shL}{\mathscr{L}}

\newcommand{\catC}{\mathcal{C}}

\newcommand{\Ad}{\operatorname{Ad}}
\newcommand{\ad}{\operatorname{ad}}

\newcommand{\Aff}{\text{\rm Aff}}
\newcommand{\alg}{\text{\rm alg}}

\newcommand{\Aut}{\operatorname{Aut}}

\newcommand{\Br}{\operatorname{Br}}

\newcommand{\Der}{\operatorname{Der}}

\newcommand{\edim}{\operatorname{edim}}
\newcommand{\End}{\operatorname{End}}

\newcommand{\Frac}{\operatorname{Frac}}
\newcommand{\frob}{{(p)}}

\newcommand{\GL}{\operatorname{GL}}
\newcommand{\liegl}{\operatorname{\mathfrak{gl}}}

\newcommand{\Hom}{\operatorname{Hom}}
\newcommand{\Hilb}{\operatorname{Hilb}}

\newcommand{\id}{{\operatorname{id}}}

\newcommand{\Kernel}{\operatorname{Ker}}

\newcommand{\kod}{\operatorname{kod}}

\newcommand{\Lie}{\operatorname{Lie}}

\newcommand{\lra}{\longrightarrow}

\newcommand{\Mat}{\operatorname{Mat}}

\newcommand{\maxid}{\mathfrak{m}}

\renewcommand{\O}{\mathscr{O}}

\newcommand{\op}{\text{\rm op}}

\newcommand{\perf}{{\text{\rm perf}}}

\newcommand{\PGL}{\operatorname{PGL}}

\newcommand{\pmap}{{[p]}}
\newcommand{\pr}{\operatorname{pr}}

\newcommand{\quadand}{\quad\text{and}\quad}

\newcommand{\ra}{\rightarrow}

\newcommand{\Rad}{\operatorname{Rad}}
\newcommand{\rank}{\operatorname{rank}}
\newcommand{\red}{{\operatorname{red}}}

\newcommand{\sep}{{\operatorname{sep}}}
\newcommand{\Set}{{\text{\rm Set}}}
\newcommand{\Sch}{\text{\rm Sch}}

\newcommand{\Sing}{\operatorname{Sing}}
\newcommand{\SL}{\operatorname{SL}}
\newcommand{\liesl}{\operatorname{\mathfrak{sl}}}

\newcommand{\Spec}{\operatorname{Spec}}

\newcommand{\SO}{\operatorname{SO}}

\newcommand{\triv}{\text{{\rm triv}}}

\newcommand{\uHom}{\underline{\operatorname{Hom}}}

\newcommand{\lieg}{\mathfrak{g}}
\newcommand{\liea}{\mathfrak{a}}
\newcommand{\lieb}{\mathfrak{b}}
\newcommand{\lieh}{\mathfrak{h}} 
\newcommand{\lien}{\mathfrak{n}}

\newcommand{\trv}{\text{\rm triv}}

\newcommand{\inert}{\text{\rm inert}}

\begin{document}

\title[Frobenius kernels for automorphism group schemes]
      {The structure of Frobenius kernels for automorphism group schemes}

\author[Stefan Schr\"oer]{Stefan Schr\"oer}
\address{Mathematisches Institut, Heinrich-Heine-Universit\"at,
40204 D\"usseldorf, Germany}
\curraddr{}
\email{schroeer@math.uni-duesseldorf.de}

\author[]{Nikolaos Tziolas}
\address{Department of Mathematics and Statistics,
University of Cyprus,
P.O.\ Box 20537, Nicosia, Cyprus}
\curraddr{}
\email{tziolas@ucy.ac.cy}


\dedicatory{Second revised version, 15 June 2022}

\begin{abstract}
We establish   structure results for   Frobenius kernels 
of   automorphism group schemes for surfaces of general type 
in positive characteristics. It turns out that there are surprisingly few possibilities.
This relies on properties of the famous Witt algebra, which is a simple Lie algebra without finite-dimensional
counterpart over the complex numbers, together with is twisted forms.
The result actually holds true for arbitrary proper integral schemes 
under the assumption that the Frobenius kernel has large
isotropy group at the generic point. This property is measured by a new numerical invariant
called the foliation rank.
\end{abstract}

\maketitle
\tableofcontents

\section*{Introduction}
\mylabel{Introduction}

Let $k$ be an algebraically closed  ground field of characteristic $p\geq 0$ and $X$ be a proper   scheme. Then the \emph{automorphism group scheme}
$\Aut_{X/k}$ is locally of finite type, and the connected component $\Aut_{X/k}^0$ 
is of finite type. The corresponding 
Lie algebra $\lieh=H^0(X,\Theta_{X/k})$ is the   space of global vector fields.
If $X$ is smooth and of general type, then the group $\Aut(X)$ is actually finite,    according to 
a general result of Martin-Deschamps and Lewin-M\'en\'egaux
\cite{Martin-Deschamps; Lewin-Menegaux 1978}. 

Throughout this paper we are mainly interested in  characteristic $p>0$. 
Then the group scheme $\Aut_{X/k}$ comes with a relative Frobenius map,
and the resulting Frobenius kernel $H=\Aut_{X/k}[F]$ is a \emph{height-one  group scheme}.
The group  of rational points is trivial, but the coordinate ring may contain nilpotent elements. 
The Lie algebra $\lieh=H^0(X,\Theta_{X/k})$ remains the   space of  global vector fields,
or equivalently the space of $k$-linear derivations $D:\O_X\ra \O_X$. The $p$-fold composition 
in the associative ring of $k$-linear differential operators 
endows the Lie algebra with an additional structure, the so-called \emph{$p$-map} $D\mapsto D^\pmap$, which
turns $\lieh$ into a \emph{restricted Lie algebra}.
By the \emph{Demazure--Gabriel Correspondence},   height-one group schemes and restricted Lie algebras
determine each other.

Our goal is to  uncover   the   \emph{structural properties} of the height-one group scheme $H=\Aut_{X/k}[F]$,
or equivalently the restricted Lie algebra $\lieh=H^0(X,\Theta_{X/k})$, and our initial motivation
was to understand the case of surfaces of general type.
Such surfaces with    $\lieh\neq 0$ were first constructed
Russell \cite{Russell 1984} and Lang \cite{Lang 1983}. These constructions rely on Tango curves  
\cite{Tango 1972}, and come with a purely inseparable
covering by a ruled surface.
By a similar construction with abelian surfaces,
Shepherd-Barron produced examples in characteristic $p=2$ that are non-uniruled (\cite{Shepherd-Barron 1996}, Theorem 5.3). 
Ekedahl  already had  examples    with rational double points  for arbitrary $p>0$ (\cite{Ekedahl 1987}, pages 145--146);
the vector fields, however,  do not extend to a resolution of singularities.
Recently, Martin studied infinitesimal automorphism group schemes of elliptic and quasielliptic surfaces
(\cite{Martin 2020a}, \cite{Martin 2020b}).

However,  almost nothing seems to be known about the general structure of the height-one group schemes $H=\Aut_{X/k}[F]$, and one would 
expect little restrictions in this respect. The main result of this paper asserts that under certain assumptions,   quite the opposite is true:

\begin{maintheorem}
(See Thm.\ \ref{structure result})
Let $X$ be a proper integral scheme with foliation rank $r\leq 1$.  Then the Frobenius kernel $H=\Aut_{X/k}[F]$ is isomorphic to the Frobenius
kernel  of one of the following three basic types of group schemes:
$$
\SL_2\quadand \GG_a^{\oplus n} \quadand \GG_a^{\oplus n}\rtimes\GG_m,
$$
for some integer $n\geq 0$.
\end{maintheorem}

In the latter two cases, the respective Frobenius kernels are $\alpha_p^{\oplus n}$ and the semidirect product $\alpha_p^{\oplus n}\rtimes\mu_p$.
The \emph{foliation rank} is a new invariant  that can  be defined as follows:
Forming the quotient $Y=X/H$ by the Frobenius kernel of the automorphism group scheme, 
the canonical map $X\ra Y$ induces a height-one extension $E=k(Y)\subset k(X)=F$
of function fields, and the foliation rank $r\geq 0$ is given by $[F:E]=p^r$. Via the  Jacobson Correspondence,
this can also be   expressed in terms of the  \emph{inertia subgroup scheme} for the induced action of the base-change $H_F$
on $F\otimes_EF$. This geometric interpretation of the Jacobson Correspondence seems to be  of independent interest (see Section \ref{Inertia}). 

If $X$ is a proper normal surface with $h^0(\omega_X^\vee)=0$, for example a surface of general type or
a properly elliptic surface, the foliation rank is automatically $r\leq 1$, 
and the above result applies  (see Corollary \ref{structure for surfaces}).  Indeed, our initial motivation was to find restrictions on the Frobenius kernels
for surfaces of general type.

The key idea in the proof  is  to relate our geometric   problem to algebraic properties of  the famous \emph{Witt algebra} $\lieg_0=\Der_E(F_0)$,
formed with  the truncated polynomial ring $F_0=E[t]/(t^p)$ over certain function fields $E$. 
This algebra was indeed introduced by Ernst Witt, compare
the discussion in \cite{Strade 1993}.  It is one of the simple algebras  in odd characteristic $p>0$   having no
finite-dimensional counterpart over the complex numbers.
Note that it has nothing to do with  the ring of Witt vectors, or Witt groups for quadratic forms.

The foliation rank is $r=1$ if and only if  $\deg(X/Y)=p$. This situation is  
paradoxical, because it may hold even with large Frobenius kernels.
We now   compare $H= \Aut_{X/k}[F]$ with the generic fiber
of the relative group scheme $G=\Aut_{X/Y}$. In other words, we relate the restricted Lie algebra 
$\lieh=H^0(X,\Theta_{X/k})$ over $k$ with the restricted Lie algebra $\lieg=\Der_E(F)$ over the function field
$E=k(Y)$.
The latter is a \emph{twisted form} of the  Witt algebra  $\lieg_0=\Der_E(F_0)$.
The classification of its subalgebras  due to Premet and
Stewart \cite{Premet; Stewart 2019} is one key ingredient for our proof.
Among other surprising  features, $\lieg_0$ contains Cartan algebras of different dimensions.
A crucial  observations is that \emph{the bigger Cartan algebras disappear after passing to twisted forms
like $\lieg$, leaving few possibilities for subalgebras}.
This is an algebraic incarnation for the fact that the reduced part of a group scheme may not
be a subgroup scheme, and if it is, it may not be normal.

The semidirect products $\alpha_p^{\oplus n}\rtimes\mu_p$ indeed occur as Frobenius kernels of automorphism group schemes.
In Section \ref{Examples}, we construct examples of surfaces   as  coverings $X\ra \PP^2$ of degree $p$ or divisors $X\subset\PP^3$ of degree $2p+1$,
such that $\lieh=k^n\rtimes\liegl_1(k)$, for certain integers $n\geq 0$. So far, we do not know if $\lieh=\liesl_2(k)$ may also occur. 
In our examples, the minimal resolutions are surfaces $S$ of general type, and $X$ are their canonical models.

Such $X$ are also called \emph{canonically polarized surface}.
They  come  with two \emph{Chern numbers} $c^2_1=c_1^2(L_{X/k}^\bullet)=K_X^2$
and $c_2=c_2(L_{X/k}^\bullet)$. This was introduced by Ekedahl, Hyland and Shepherd-Barron
\cite{Ekedahl; Hyland; Shepherd-Barron 2012} for general  proper surfaces whose local rings are complete intersections, such
that the cotangent complex   is perfect.
Using Noether's inequality and results from Ekedahl \cite{Ekedahl 1988}, we show with more classical   methods:

\begin{maintheorem}
(see Thm.\ \ref{bound result})
Let $X$ be canonically polarized surface, with Chern numbers $c_1^2$ and $c_2$.
Then  the Lie algebra $\lieh=H^0(X,\Theta_{X/k})$ for the  Frobenius kernel $H=\Aut_{X/k}[F]$ has the property
$\dim(\lieh)\leq \Phi(c_1^2,c_2)$ for the polynomial
$$
\Phi(x,y) =
\begin{cases}
\frac{1}{144}(73x+y)^2-1, & \mathrm{if} \; c_1^2\geq 2\\
\frac{1}{144}(121x+y)^2-1, & \mathrm{if}\; c_1^2=1
\end{cases}
$$
\end{maintheorem}
With Noether's Inequality, this also gives the weaker bound $\dim(\lieh)\leq \Psi(c_1^2)$ with the polynomial 
$$
\Psi(x)=
\begin{cases}
\frac{169}{4}x^2+39x+8		& \text{if $c_1^2\geq 2$;}\\
\frac{441}{4}x^2 + 63x+8	& \text{if $c_1^2=1$.}
\end{cases}
$$
Note that Xiao  \cite{Xiao 1995} proved  $|\Aut(X)|\leq 1764 c_1^2$ over the complex numbers.

\medskip
The paper is organized as follows:
Section \ref{Restricted} contains  general facts on restricted Lie algebras
and their semidirect products. In Section \ref{Toral rank} we examine
multiplicative and additive vectors, and  the toral rank.
In Section \ref{Automorphism} we collect general facts on automorphism group schemes
for proper schemes,  the quotient 
by height-one group schemes, and discuss   twisted forms of some relevant restricted Lie algebras.
Section \ref{Inertia} contains a geometric interpretation of the Jacobson
correspondence, in terms of inertia group schemes at generic points.
We introduce the foliation rank and establish its basic properties in Section \ref{Foliation rank}.
In Section \ref{Invariant subspaces} we analyze the removal of subvector spaces
under certain   twists. Then we make a detailed analysis of the automorphism group scheme
for radical extensions of prime degree in Section \ref{Radical}, followed by an examination
of the corresponding Witt algebras in Section \ref{Witt algebras}.
In Section \ref{Twisting} we show how structural properties of restricted Lie algebras over
different fields are inherited. Our main result on the structure of the Frobenius kernel
for automorphism groups is contained in Section \ref{Subalgebras}.
 Section \ref{Canonically polarized} contains the bound for surfaces of general type.
In the final Section \ref{Examples}, we construct examples.

\begin{acknowledgement}
The research started during a visit of the first author at the University of Cyprus, and he
likes to thank the Department of Mathematics and Statistics for its hospitality. 
The research was also conducted       in the framework of the   research training group
\emph{GRK 2240: Algebro-Geometric Methods in Algebra, Arithmetic and Topology}.
We like
to thank Yuya Matsumoto and Michel Brion for   valuable comments. We also thank the referees for
many helpful remarks, and in particular for  the suggestion  to use Tango curves in Section \ref{Examples}.

\end{acknowledgement}

\section{Restricted Lie algebras}
\mylabel{Restricted}

In this section we review   some standard results on restricted Lie algebras and height-one group schemes that
are relevant for the applications we have in mind.
Let $k$ be a ground field of characteristic $p>0$.
For each   ring $R$, not necessarily commutative or associative, the vector space  $\Der_k(R)$ of $k$-derivations $D:R\ra R$
is closed under forming commutators $[D,D']$ and $p$-fold compositions $D^p$
in the associative ring $\End_k(R)$.
One now views   $\Der_k(R)$ as a \emph{Lie algebra}, endowed the map $D\mapsto D^p$ as an additional structure.

This leads to the following abstraction:
A \emph{restricted Lie algebra} is a Lie algebra $\lieg$, together with a map $\lieg\ra \lieg$, $x\mapsto x^\pmap$
called the \emph{$p$-map}, subject to the following three axioms: 
\begin{enumerate}
\item[(R 1)]
We have $\ad_{x^\pmap} = (\ad_x)^p$ for all vectors $x\in \lieg$.
\smallskip
\item[(R 2)]
Moreover $(\lambda\cdot x)^\pmap = \lambda^p\cdot x^\pmap$ for all  vectors $x\in \lieg$ and scalars $\lambda\in k$.
\smallskip
\item [(R 3)]
The formula $(x+y)^\pmap = x^\pmap + y^\pmap + \sum_{r=1}^{p-1}s_r(x,y)$ holds for all $x,y\in\lieg$.
\end{enumerate}
Here the summands $s_r(x,y)$ are   universal expressions defined by 
$$
s_r(t_0,t_1)=-\frac{1}{r}\sum_u  (\ad_{t_{u(1)}}\circ \ad_{t_{u(2)}}\circ \ldots \circ\ad_{t_{u(p-1)}} )(t_1),
$$
where $\ad_a(x)=[a,x]$ denotes the \emph{adjoint representation}, and the 
the index runs over all maps $u:\{1,\ldots,p-1\}\ra \{0,1\}$ taking the value zero exactly $r$ times.
For $p=2$ the expression simplifies to $s_1=[t_0,t_1]$, whereas $p=3$ gives $s_1=[t_1,[t_0,t_1]]$ and $s_2=[t_0,[t_0,t_1]]$.
Restricted Lie algebras were  introduced and studied by Jacobson \cite{Jacobson 1937},   and   also  go  under the name   \emph{$p$-Lie algebras}.
We refer to the monographs of Demazure and Gabriel \cite{Demazure; Gabriel 1970}, in particular Chapter II, \S7,
or Strade and Farnsteiner \cite{Strade; Farnsteiner  1988} for more details.

Throughout the paper, terms like homomorphisms, subalgebras, ideals, extensions etc.\ are understood in the \emph{restricted  sense}, if not said
otherwise. For example, an \emph{ideal} $\ideala\subset\lieg$ is a  vector subspace such that $[x,y], x^\pmap\in\ideala$ whenever
$x\in \ideala$ and $y\in \lieg$. Note that this holds    for  the \emph{center} $\mathfrak{C}(\lieg)=\{a\in\lieg\mid \text{$[a,x]=0$ for all $x\in \lieg$}\}$,
because $[a^\pmap,x]=(\ad_a)^p(x) =(\ad_a)^{p-1}([a,x])=0$. 

For abelian $\lieg$,   the $p$-map becomes \emph{semi-linear}, which means that it corresponds to a linear map  $\lieg\ra\lieg$
when the scalar multiplication in the range is redefined via    Frobenius.  In turn, those $\lieg$
correspond to modules over the associative polynomial ring $k[F]$, in which the relation $F\lambda=\lambda^pF$ holds.
Every right ideal is principal; this also holds for left ideals, provided that $k$ is perfect, and then the structure theory 
developed by  Jacobson applies (\cite{Jacobson 1943}, Chapter 3).

In contrast, for non-abelian $\lieg$ 
the $p$-map \emph{fails to be additive}, and it is   challenging to understand its structure.
However, by axiom (R 1) it is determined by the bracket  up to central elements, because
$ [a^\pmap,x] = (\ad_a)^p(x)$. In particular, if the center is trivial,
the $p$-map is unique, once it exists. This also explains the terminology \emph{restricted}.
 
Recall that for each group scheme $G$, the \emph{Lie algebra} $\lieg=\Lie(G)$ is defined by the short exact sequence
$$
0\lra \Lie(G)\lra G(k[\epsilon])\lra G(k)\lra 0,
$$
where $k[\epsilon]$ is the ring of dual numbers, and the map is the restriction with respect to the inclusion $k\subset k[\epsilon]$.
As explained in \cite{Demazure; Gabriel 1970}, Chapter II, \S7, it carries the structure of a restricted Lie algebra, in a functorial way.
Also recall  that the relative Frobenius   $F:G\ra G^{(p)}$ is a homomorphism. The resulting \emph{Frobenius kernel} $G[F]$
is a group scheme whose underlying topological space is a singleton.

Let us call $G$   \emph{of height one} if it is of finite type  and annihilated by the relative Frobenius map.
We then also say that $G$ is a \emph{height-one group scheme}.
According to \cite{Demazure; Gabriel 1970}, Chapter II, \S7, Theorem 3.5 the canonical map 
$$
\Hom(G,H)\lra \Hom(\Lie(G),\Lie(H))
$$
is bijective whenever $G$ has height   one.  In particular,
the functor $G\mapsto \Lie(G)$ is an equivalence   between the category of height-one group schemes 
and the category of finite-dimensional restricted Lie algebras. We call this the \emph{Demazure--Gabriel Correspondence}.
The inverse functor sends $\lieg$ 
to the spectrum of the dual for the Hopf algebra $U^\pmap(\lieg)$, which is the
universal enveloping algebra $U(\lieg)$ modulo the ideal generated by the elements $x^p-x^\pmap$, for $x\in\lieg$.
From this one deduces the formulas
$$
|G|=h^0(\O_G) = p^{\dim(\lieg)}\quadand \edim(\O_{G,e})=\dim(\lieg).
$$
As customary, we write $\liegl_n(k)$ for the restricted Lie algebra of $n\times n$-matrices, where bracket and $p$-map 
are given by commutators and $p$-powers,  and  $\liesl_n(k)$ for the ideal of trace zero matrices. Furthermore, $k^n$ denotes
the standard vector space, endowed with trivial bracket and $p$-map. 

Let $\ideala\subset\lieg$ be an ideal, and consider the vector space $\Der_k(\liea) $   of all $k$-linear
derivations. Then $\Der_k(\liea)\subset\liegl(\liea)$ is a subalgebra.
Derivations $D:\lieg\ra \lieg$ satisfying the additional condition $D(a^\pmap)=(\ad_a)^{p-1}(D(a))$ for all $a\in \liea$ are called
\emph{restricted derivations}. Write $\Der'_k(\liea)$ for the vector space of all restricted $k$-derivations. According to  \cite{Jacobson 1941b}, Theorem 4
the inclusion $\Der'_k(\liea)\subset\Der_k(\liea)$ is a  subalgebra. 
By the Jacobi identity and axiom (R 1), the adjoint map   defines a homomorphism
\begin{equation}
\label{restriced homomorphism}
\lieg\lra \Der'_k(\liea),\quad x\longmapsto (a\mapsto [x,a]).
\end{equation}
Given restricted Lie algebras $\lieh$ and $\liea$, we are  are now interested in \emph{extensions}
$$
0\lra \liea\lra\lieg\lra \lieh\lra 0,
$$
such that  $\liea\subset \lieg$ becomes  an ideal with quotient $\lieg/\liea=\lieh$.
The   extension \emph{splits} if the ideal $\ideala\subset\lieg$
admits a complementary subalgebra $\lieh'\subset\lieg$. 
Composing the inverse for the projection   $\lieh'\ra \lieh$ with 
\eqref{restriced homomorphism}, we obtain  a homomorphism $\varphi:\lieh\ra\Der'_k(\liea)$.
Conversely, suppose we have such a homomorphism, written as $h\mapsto(a\mapsto\varphi_h(a))$.
On the vector space sum $\ideala\oplus\lieh$, we now define   bracket and $p$-map by  
\begin{equation}
\begin{gathered}
{}
[a+h,a'+h'] = [a,a'] + [h,h'] +  \varphi_h(a') - \varphi_{h'}(a),\\
(a+h)^\pmap = a^\pmap + h^\pmap + \sum_{r=1}^{p-1}s_r(a,h).
\end{gathered}\label{structure semidirect}
\end{equation}

\begin{lemma}
\mylabel{split extensions}
The  above  endows the vector space $\lieg=\ideala\oplus\lieh$ with the structure of a restricted Lie algebra,
such that $\ideala$ and $\idealh$ is an ideal and subalgebra, respectively.
\end{lemma}

\proof
As explained in  \cite{LIE 1-3}, Chapter I, \S1.8  the bracket turns
$\lieg=\ideala\oplus\lieh$ into a Lie algebra, having $\ideala$ as an ideal and $\idealh$ as a subalgebra.
Now choose   bases $a_i\in\ideala$ and $h_j\in \idealh$, such that   $a_i,h_j$ form a basis for $\lieg$.
We claim that 
\begin{equation}
\label{jacobson condition}
(\ad_{a_i})^p=\ad_{(a_i)^\pmap}\quadand (\ad_{h_j})^p=\ad_{(h_j)^\pmap}
\end{equation}
as $k$-linear endomorphisms of $\lieg$. Indeed: 
Since $\liea$ and $\lieh$ are restricted, and by the definition of the bracket in $\lieg$, it is enough to verify
$(\ad_{a_i})^p(h)=-\varphi_h( (a_i)^\pmap)$ for every vector $h\in \lieh$, and 
$(\ad_{h_j})^p(a)=\varphi_{(h_j)^\pmap}(a)$ for every $a\in \liea$.
Since the derivations $\varphi_h$ are restricted, we   have
$$
-\varphi_h((a_i)^\pmap) = -\ad_{a_i}^{p-1}(\varphi_h(a_i)) =- \ad_{a_i}^{p-1}([h,a_i])=(\ad_{a_i})^p(h).
$$
The argument for $(\ad_{h_j})^p(a)$ is similar. Thus \eqref{jacobson condition} holds.
According to \cite{ Strade; Farnsteiner  1988}, Theorem 2.3 there is a unique $p$-map satisfying
\eqref{jacobson condition} and the   axioms (R 1)--(R 3). 
By construction, this $p$-map on $\lieg$ coincides with the given $p$-map on $\liea$ and $\lieh$.
It thus coincides with \eqref{structure semidirect}, in light of the third axiom.
\qed

\medskip
In the above situation,  the restricted Lie algebras  $\lieg=\liea\rtimes_\varphi\lieh$ are called 
\emph{semidirect products}. Obviously,  every split extension of $\lieh$ by $\liea$ is of this form.
Of particular importance for us is the case $\liea=k^n$ and $\lieb=\liegl_1(k)$, where the homomorphism $\varphi:\liegl_1(k)\ra\liegl(k^n)=\Der_k(k^n)$
sends scalars to scalar matrices. The resulting restricted Lie algebra is written as $k^n\rtimes\liegl_1(k)$.
Here bracket and $p$-map are given by the formulas 
\begin{equation}
\label{semidirect special}
[v+\lambda e,v'+\lambda'e]=\lambda v'-\lambda'v\quadand (v+\lambda e)^\pmap = \lambda^{p-1}(v+\lambda e),
\end{equation}
where $e\in\liegl_1(k)$ is the unit element,  and  $v,v'\in k^n$ are vectors, and $\lambda,\lambda'\in k$ are scalars.

\section{Toral rank and \texorpdfstring{$p$}{p}-closed vectors}
\mylabel{Toral rank}

Let $\lieg$ be a finite-dimensional restricted Lie algebra over a ground field $k$ of characteristic $p>0$, 
and $G$ be the corresponding height-one group scheme,
such that $\Lie(G)=\lieg$.
Recall that   $x\in \lieg$ is called \emph{$p$-closed} if $x^{[p]}\in kx$. 
Such vectors are   called \emph{multiplicative} if $x^\pmap \neq 0$, and \emph{additive} if $x^\pmap =0$.
If the vector is non-zero,   $\lieh=kx$ is a one-dimensional subalgebra, hence
corresponds to a subgroup scheme $H\subset G$ of order $p$.
For multiplicative vectors, this is a twisted form of the diagonalizable group scheme $\mu_p=\GG_m[F]$.
In the additive case, it is isomorphic to the unipotent group scheme $\alpha_p=\GG_a[F]$.
This basic fact has many geometric applications: For results concerning K3 surfaces, Enriques surfaces and  Kummer surfaces, see 
\cite{Schroeer 2007}, \cite{Schroeer 2017} and \cite{Kondo; Schroeer 2019}.

\begin{proposition}
\mylabel{some p-closed}
Every vector in $\lieg=k^n\rtimes\liegl_1(k)$ is $p$-closed.
The same holds for $\lieg=\liesl_2(k)$ in characteristic $p\geq 3$. 
\end{proposition}

\proof
The first assertion immediately follows from \eqref{semidirect special}.  
Recall that $\liesl_2(k)$ is the restricted Lie algebra comprising the traceless matrices 
$
A=(\begin{smallmatrix}a&b\\c&-a\end{smallmatrix})\in\Mat_2(k)$.
The characteristic polynomial   $\chi_A(T)=T^2+d$ depends only on the  determinant $d=-a^2-bc$,
so the possible Jordan normal forms over $k^\alg$ are
$$
\begin{pmatrix} \sqrt{d}&0\\0&-\sqrt{d}\end{pmatrix}\quadand \begin{pmatrix}0&0\\1&0 \end{pmatrix}.
$$
Computing $p$-powers via the above normal forms, we see that  $A^\pmap= d^{(p-1)/2}A$.
\qed

\medskip
The traceless matrices 
$h=(\begin{smallmatrix}1&0\\0&-1\end{smallmatrix})$ and $x=(\begin{smallmatrix}0&1\\0&1\end{smallmatrix})$ and 
$y=(\begin{smallmatrix}0&0\\1&0\end{smallmatrix})$ form a basis of $\liesl_2(k)$, and the structural constants 
are given by
$$
[h,x]=2x,\quad [h,y]=2y,\quad [x,y]= h^\pmap=h,\quad x^\pmap=y^\pmap =0.
$$
One also says that $(h,x,y)$ is an \emph{$\liesl_2(k)$-triple}. For $p\geq 3$, it follows that for each non-zero $a\in\liesl_2(k)$,
the adjoint map $\ad_a$ is bijective, hence $\liesl_2(k)$ is simple.
In contrast, for   $p=2$ we have a central extension $0\ra \liegl_1(k)\ra\liesl_2(k)\ra k^2\ra 0$,
where the kernel corresponds to scalar matrices.  The
extension does not split, because $A^{[2]}\neq 0$ for all matrices not contained in the kernel.

If $k$ is algebraically closed, the   \emph{toral rank}  for a restricted Lie algebra $\lieg$ 
is the maximal integer $r\geq 0$ for which there is an embedding $\liegl_1(k)^{\oplus r}\subset \lieg$.
In terms of vectors, the condition  means that there are linearly independent $x_1,\ldots,x_r\in \lieg$ with $[x_i,x_j]=0$ and $x_i^\pmap=x_i$.
For general fields $k$, we define the toral rank  as the toral rank of the base-change $\lieg\otimes_kk^\alg$.
Following the notation in  \cite{SGA 3b}, Expos\'e XII, Section 2 we denote  this integer by $\rho_t(\lieg)\geq 0$. 
By Hilbert's Nullstellensatz the toral rank does not change under field extensions.
According to \cite{Block; Wilson 1988}, Lemma 1.7.2  it satisfies $\rho_t(\lieg)=\rho_t(\lien)+\rho_t(\lieg/\lien)$ 
for each ideal $\ideala\subset\lieg$. In other words, it is additive in extensions.
Obviously $0\leq \rho_t(\lieg)\leq \dim(\lieg)$. 

\begin{proposition}
The following are equivalent:
\begin{enumerate}
\item The restricted Lie algebra $\lieg$ has maximal toral rank $\rho_t(\lieg)=\dim(\lieg)$.
\item The group scheme $G$ is a twisted form of some $\mu_p^{\oplus r}$.
\item The group scheme $G$ is multiplicative.
\end{enumerate}
\end{proposition}

\proof
It suffices to treat the case that $k$ is algebraically closed. The implications (i)$\Leftrightarrow$(ii)$\Rightarrow$(iii)
are obvious. Now suppose that (iii) holds. Since $k=k^\alg$ the group scheme $G$ is diagonalizable, whence the spectrum
of the Hopf algebra $k[\Lambda]$ for some finitely generated abelian group $\Lambda$. We have $p\Lambda=0$ because 
 $G$ has height  one. Choosing an $\FF_p$-basis for $\Lambda$ gives $G=\mu_p^{\oplus r}$, thus (ii) holds.
\qed

\medskip 
The other extreme is somewhat more involved:

\begin{proposition}
The following are equivalent:
\begin{enumerate}
\item The restricted Lie algebra $\lieg$ has  minimal toral rank $\rho_t(\lieg)=0$.
\item There is some exponent $\nu\geq 0$ with $x^{[p^\nu]}=0$ for all vectors $x\in\lieg$.
\item There are ideals $0=\ideala_0\subset\ldots\subset\ideala_r=\lieg$ inside $\lieg$ with quotients
$\liea_i/\liea_{i-1}\simeq k$.
\item There are normal subgroup schemes $0=N_0\subset \ldots\subset N_r=G$ inside $G$ with quotients $N_i/N_{i-1}\simeq\alpha_p$.
\item The group scheme $G$ is unipotent.
\end{enumerate}
\end{proposition}

\proof

The implications (iv)$\Rightarrow$(v) and (iii)$\Rightarrow$(ii)  and (ii)$\Rightarrow$(i) are trivial, whereas
(iv)$\Leftrightarrow$(iii) follows  from the Demazure--Gabriel Correspondence.

We next verify (v)$\Rightarrow$(i). Without loss of generality we may assume that $k$ is algebraically closed.
Then there is a composition series $G_j$ inside $G$ such that $G_j/G_{j-1}$ is isomorphic to  a subgroup scheme
of the additive group $\GG_a$. This already lies in $\alpha_p=\GG_a[F]$, because   $G$ has height one.
For the corresponding subalgebras $\lieb_j$ inside $\lieg$ this means $\lieb_j/\lieb_{j-1}\subset k$.
The additivity of toral rank implies $\rho_t(\lieg)=0$.

To see (i)$\Rightarrow$(ii) we   may assume that $k$ is algebraically closed, and then the implication
follows from   \cite{Premet 1990}, Corollary 2.
For (ii)$\Rightarrow$(iii) we use $(\ad_a)^{p^\nu} =\ad_{a^{[p^\nu]}}=0$, and conclude
with Engel's Theorem (\cite{LIE 1-3}, Chapter I,   \S4.2) that the underlying Lie algebra $\lieg$ is nilpotent.
Now recall that the center $\mathfrak{C}(\lieg) $
is invariant under the $p$-map. In turn,
the upper central series, which is recursively defined by $\lieg_{i+1}/\lieg_i = \mathfrak{C}(\lieg/\lieg_i)$,
yields  a sequence of ideals   $0=\lieg_0\subset\ldots\subset\lieg_s=\lieg$ having abelian quotients.
This reduces our problem to the case that $\lieg$ itself is abelian.
We   proceed by induction on $n=\dim(\lieg)$. The case $n=0$ is trivial. Suppose now that $n>0$, 
and that (iii) holds for $n-1$. Fix some   $x\neq 0$, and consider the largest exponent $d\geq 1$ such that  $x^{[p^d]}\neq 0$.
Replacing $x$ by $x^{[p^d]}$, we may assume that $x^\pmap =0$. Then $\liea_1=kx$ is a one-dimensional ideal.
The quotient $\lieg'=\lieg/\ideala$ has dimension $n'=n-1$, and furthermore $\rho_t(\lieg')=0$ by additivity of toral rank.
To the induction hypothesis applies to $\lieg'=\lieg/\liea$, and the Isomorphism Theorem
gives the desired ideals in $\lieg$. 
\qed

\section{Automorphism group schemes}
\mylabel{Automorphism}

Let $k$ be a ground field. Write $(\Aff/k)$ for the category of affine $k$-schemes, which we usually write as $T=\Spec(R)$.
Recall that an \emph{algebraic space} is a contravariant functor $X:(\Aff/k)\ra(\Set)$
satisfying the sheaf axiom with respect to the \'etale topology, such that   
the diagonal $X\ra X\times X$ is relatively representable by schemes,
and that there is an \'etale surjection $U\ra X$ from some scheme $U$. 
According to \cite{SP}, Lemma 076M, the sheaf axiom already holds with respect to the fppf topology.
Throughout, we use the fppf topology if not stated otherwise.
Algebraic spaces  are important
generalizations of schemes, because modifications, quotients, families, or moduli spaces  of schemes  
are frequently algebraic spaces rather than schemes.  
We refer to the monographs of 
Olsson \cite{Olsson 2016}, Laumon and Moret-Bailly \cite{Laumon; Moret-Bailly 2000}, 
Artin \cite{Artin 1971}, Knutson \cite{Knutson 1971}, and to the stacks project
\cite{SP}, Part 4. 

Let $X$ be a scheme, or more generally an algebraic space, that is separated and of finite type.
Recall that the $R$-valued points of the  Hilbert functor $\Hilb_{X/k}$
are the closed subschemes $Z\subset X\otimes R$ such that the projection $Z\ra\Spec(R)$ is proper and flat.
Regarding automorphisms $f:X\otimes R\ra X\otimes R$ as graphs, we see that $\Aut_{X/k}$ is an open subfunctor.
According to \cite{Artin 1969}, Theorem 6.1 the Hilbert  functor  is representable by an algebraic space that  is separated and locally of finite type.
In turn, the same holds for $\Aut_{X/k}$, which additionally carries 
a group structure. Using descent and translations, one sees that it must be schematic.
 The  Lie algebra for the automorphism group scheme  is given by 
$$
\Lie(\Aut_{X/k})= H^0(X,\Theta_{X/k}),
$$
where $\Theta_{X/k}=\uHom(\Omega^1_{X/k},\O_X)$ is the coherent sheaf dual to the sheaf of K\"ahler differentials.
 
We now assume that $X$ is proper, and that the ground field has characteristic $p>0$. 
Then $\lieg=H^0(X,\Theta_X)$  is a restricted Lie algebra of finite dimension,
which corresponds to the Frobenius kernel $G[F]$ for the automorphism group scheme $G=\Aut_{X/k}$.
Note that $G[F]$ is a height-one group scheme, of order $p^n$, where $n=h^0(\Theta_{X/k})$.

Let $H$ be a group scheme that is separated and locally of finite type, $f:H\ra G$ be a homomorphism, and $P$ be a $H$-torsor.
The latter is an algebraic space, endowed with a free and transitive $H$-action. The set of isomorphism classes comprise
the non-abelian cohomology $H^1(k,H)$, formed with respect to the fppf topology.
On the product $P\times X$ we get a diagonal action. This action is free, because it is free on the first factor.
It follows that the quotient ${}^P\!X=H\backslash(P\times X)$ exists as an algebraic space (see for example \cite{Laurent; Schroeer 2021}, Lemma 1.1).
We have ${}^P\!X\simeq X$ provided that $P$ is trivial, that is, contains a rational point.
In any case, there is an \'etale surjection $U\ra P$ from some scheme $U$. According to Hilbert's Nullstellensatz,
every closed point $a\in U$ defines a finite field extension $k'=\kappa(a)$,
and we see that ${}^P\!X\otimes k'\simeq X\otimes k'$. We therefore say that  ${}^P\!X$ is a \emph{twisted form} of $X$.
Indeed, every algebraic space  $Y$ that becomes isomorphic to $X$ after some field extension is of this form,
with $H=\Aut_{X/k}$.
 
Our $f:H\ra G=\Aut_{X/k}$ induces a homomorphism $c:H\ra\Aut_{G/k}$, which sends $h\in H(R)$
to the inner automorphism $g\mapsto f(h)  g f(h)^{-1}$. This gives a twisted form ${}^P\!G$ of $G$,
and its Lie algebra ${}^P\!\lieg$  is a twisted form of $\lieg$.
In fact, one may view $\lieg$ as a vector scheme as in Section \ref{Invariant subspaces}, regard bracket and $p$-map as morphisms of  schemes, 
and obtains ${}^P\!\lieg$ by taking the rational points on the twisted form of the vector scheme, formed   via the derivative
$c':H\ra \Aut_{\lieg/k}$.

\begin{lemma}
\mylabel{twisted automorphisms}
There is a canonical identification
${}^P\!\Aut_{X/k}=\Aut_{{}^P\!X/k} $,
where on the left we take twist with respect to $c:H\ra\Aut_{G/k}$.
The restricted Lie algebra for this group scheme is ${}^P\!\lieg$,
where we twist with respect to $c':H\ra\Aut_{\lieg/k}$.
\end{lemma}

\proof
This follows from very general considerations in \cite{Giraud 1971}, Chapter III, which can be made explicit as follows:
Consider the canonical morphism
$$
P\times \Aut_{X/k}\lra \Aut_{{}^P\!X},\quad (p,\psi)\longmapsto (H\cdot(p,x)\mapsto H\cdot(p,\psi(x)),
$$
where the description on the right is viewed as  a natural transformation for $R$-valued points.
This is well-defined, because in presence of $p\in P(R)$ the projection $\{p\}\times X(R)\ra ({}^P\!X)(R)$ is bijective.
For each $h\in H(R)$, the element $(hp, h\psi h^{-1})$ sends the orbit $H\cdot(p,x)=H\cdot(hp,hx)$ to the orbit $H\cdot(hp,h\psi (x))=H\cdot(p,\psi(x))$.
Thus the above transformation descends to a morphism ${}^P\!\Aut_{X/k}\ra\Aut_{{}^P\!X}$, where $H$ acts via conjugacy
on $\Aut_{X/k}$. The same argument applies for the Frobenius kernel, and equivalently to the restricted Lie algebra.
\qed

\medskip
We now change notation, and suppose that $G=X$ is a height-one group scheme, and write $\lieg=\Lie(G)$.
One easily  checks that $\Aut_{G/k}$ is a closed subgroup scheme of the general linear group $\GL_{V/k}$, where
$V=H^0(G,\O_G)$.
By the Demazure--Gabriel Correspondence, used in the relative form, we get an identification $\Aut_{G/k}=\Aut_{\lieg/k}$.
The latter can be constructed directly: Choose  a basis $e_1,\ldots, e_n\in \lieg$. Then $\Aut_{\lieg/k}$ is the closed subgroup scheme
inside $\GL_{k,n}$ respecting the structural equations  
$[e_r,e_s] = \sum \lambda_{r,s,i}e_i$ and $e_r^\pmap = \sum \mu_{r,j} e_j$. 
For later use, we compute some automorphism group schemes $\Aut_{\lieg/k}$:

\begin{proposition}
\mylabel{twisted forms}
The following table lists the  automorphism group schemes and  the   resulting cohomology groups or sets for the
restricted Lie algebras $k$, $\liegl_1(k)$, $k\rtimes\liegl_1(k) $ and $\liesl_2(k)$, 
where  the last column is only valid for $p\geq 3$:
$$
\begin{array}{lllll}
\toprule
\lieg   		& k	& \liegl_1(k)		& k\rtimes\liegl_1(k) 	& \liesl_2(k)\\
\midrule
\Aut_{\lieg/k}  	& \GG_m	& \mu_{p-1}			& \GG_a\rtimes\GG_m   	& \PGL_2	\\
\midrule
H^1(k,\Aut_{\lieg/k})   & \{1\}	& k^\times/k^{\times(p-1)}	& \text{\rm singleton} 		& \text{\rm subset of $\Br(k)[2]$}\\
\bottomrule
\end{array}
$$
Here  $\Br(k)[2]$ is the kernel of multiplication by two on the Brauer group.
\end{proposition}

\proof
For the first case  $\lieg=k$ we immediately get $\Aut_{\lieg/k}=\GL_1=\GG_m$,  and Hilbert 90 gives $H^1(k,\GG_m)=\{1\}$, 
at least for the \'etale toplogy. See the discussion at the  beginning of Section \ref{Invariant subspaces} for the fppf topology.

In the second case, the  restricted Lie algebra $\lieg=\liegl_1(k)$ is generated by one element $A_1$, which gives an embedding
$\Aut_{\lieg/k}\subset\GG_m$. The structure  for $\lieg$  is given by $A_1^\pmap =A_1$. For each $k$-algebra $R$ and each invertible scalar $\lambda\in R^\times$
we thus have $\lambda^p A_1^\pmap = \lambda A_1$, and thus $\lambda^{p-1}=1$. Conversely, each such $\lambda$
gives an automorphism, hence  $\Aut_{\lieg/k}=\mu_{p-1}$.
The Kummer sequence yields $H^1(k,\Aut_{\lieg/k})=k^\times/k^{\times (p-1)}$.

The restricted Lie algebra $\lieg=k\rtimes\liegl_1(k)$ is generated inside $\liegl_2(k)$ by the matrices
$A_1=(\begin{smallmatrix}0&1\\0&0\end{smallmatrix})$ and $A_2=(\begin{smallmatrix}1&0\\0&0\end{smallmatrix})$, 
which gives an embedding $\Aut_{\lieg/k}\subset\GL_2$.
For each $R$-valued point $\varphi=(\begin{smallmatrix}a&b\\c&d\end{smallmatrix})$ from the automorphism group scheme, the condition 
$[\varphi(A_1),\varphi(A_2)]=\varphi(A_1)$
implies $c=0$ and $a=ad$. It follows $a\in R^\times$ and $d=1$, and we obtain $\Aut_{\lieg/k}\subset\GG_a\rtimes\GG_m$.
Conversely, one easily sees that each   matrix with $c=0$ and $d=1$ yields an automorphism of the restricted Lie algebra.
Now let $T$ be a torsor over $k$ with respect to $\GG_a\rtimes\GG_m$. The induced $\GG_m$-torsor has a rational point,
by Hilbert 90. Its preimage $T'\subset T$ is a torsor for $\GG_a$. Over any affine scheme, the higher cohomology of $\GG_a$ 
vanishes, so $T'$ also contains a rational point, and the torsor $T$ is trivial.

We come to the last case $\lieg=\liesl_2(k)$, which is freely generated by the matrices
$A_1=(\begin{smallmatrix}1&0\\0&-1\end{smallmatrix})$ and  $A_2=(\begin{smallmatrix}0&1\\0&0\end{smallmatrix})$ and
$A_3=(\begin{smallmatrix}0&0\\1&0\end{smallmatrix})$.  This gives an inclusion $\Aut_{\lieg/k}\subset\GL_3$.
Conjugacy $A\mapsto S AS^{-1}$ yields $\PGL_2\subset\Aut_{\lieg/k}$.
We already saw in the proof for Proposition \ref{some p-closed} that  $A^\pmap = \det(A)^{(p-1)/2}A$ for all
$A=(\begin{smallmatrix}a&b\\c&-a\end{smallmatrix})$. Moreover, $\det(A)=-a^2-bc$ defines, up to sign, the standard smooth quadratic form
on $\liesl_2(k)$ viewed as the affine  space $\AA^3$, which gives  $ \Aut_{\lieg/k}\subset O(3)$.
We have  $\PGL_2\subset\SO(3)$, because the former is connected, and this inclusion is an equality because
both are smooth and three-dimensional. This shows $\SO(3)\subset\Aut_{\lieg/k}\subset O(3)$.
From $[-A_2,-A_3]= [A_2,A_3] = A_1\neq -A_1$ we conclude that $A\mapsto -A$ is not an automorphism of $\lieg$,
so $\PGL_2=\SO(3)\subset\Aut_{\lieg/k}$ must be  an equality.

Finally, we have a central extension $0\ra\GG_m\ra\GL_2\ra\PGL_2\ra 1$, and   get maps in non-abelian cohomology
$$
H^1(k,\GL_2)\lra H^1(k,\PGL_2)\lra H^2(k,\GG_m).
$$
The term on the left is a singleton, by Hilbert 90, whereas the term on the right equals the Brauer group $\Br(k)$.
It follows that the coboundary map is injective (\cite{Giraud 1971}, Chapter IV, Proposition 4.2.8), 
and its image is  contained in the 2-torsion part of the Brauer group (\cite{GB I}, Proposition 1.4).
Thus $H^1(k,\Aut_{\lieg/k})$ is a certain subset inside the group $\Br(k)[2]$.
\qed

\medskip
Note that according to the Theorem of Merkurjev  \cite{Merkurjev 1982},
the group $\Br(k)[2]$ is generated by classes from $H^1(k,\PGL_2)$.
This set of generators, however,   is not a subgroup   in general (see \cite{Gille; Szamuely 2006}, Example 1.5.7).

\section{Quotients by height-one group schemes}
\mylabel{Quotients}

Let $k$ be a ground field of characteristic $p>0$, and $G$ a height-one group scheme, with restricted Lie algebra $\lieg$.
Suppose $X$ is a    scheme endowed with a $G$-action. Taking derivatives, we obtain a homomorphism
$\lieg\ra H^0(X,\Theta_{X/k})$
of restricted Lie algebras. According to \cite{Demazure; Gabriel 1970}, Chapter II, \S7, Proposition 3.10
any such homomorphism comes from a unique $G$-action. Note that this does not require any finiteness assumption for the scheme $X$.

We now show that such actions admit  a \emph{categorical quotient}   in the category $(\Sch/k)$
(\cite{Mumford; Fogarty; Kirwan 1993}, Definition 0.5). 
To this end we temporarily change notation  and  write   the   schemes in question
as   pairs, comprising a topological space  and a structure sheaf. Our task is to construct the categorical quotient    $(Y,\O_Y)$
for the action on $(X,\O_X)$.
First recall that the image $\O_X^p$ of the homomorphism $\O_X\ra \O_X$, $f\mapsto f^p$ is a quasicoherent $\O_X$-algebra,
with algebra structure $f\cdot g^p=(fg)^p$. In turn, the ringed space 
$(X,\O_X^p)$ is a $\FF_p$-scheme. Choose a vector space basis   $D_1,\ldots,D_n\in\lieg$.
The canonical inclusion $\O_X^p\subset\O_X$ turns $\O_X$ into a quasicoherent $\O_X^p$-algebra,
and yields the absolute Frobenius morphism $(X,\O_X)\ra(X,\O_X^p)$.
The derivations $D_i:\O_X\ra\O_X$ are $\O_{X}^p$-linear, and we write 
$\O_X^\lieg=\bigcap_{i=1}^n\Kernel(D_i)$ for the intersection of   kernels.
This is another quasicoherent $\O_X^p$-algebra.
Setting $Y=X$ and $\O_Y=\O_X^\lieg$, we obtain a scheme $(Y,\O_Y)$ that is affine over $(X,\O_X^p)$.
The identity $\id:X\ra Y$ and the canonical inclusion $\iota:\O_Y\subset\O_X$ 
define a morphism of  $\FF_p$-schemes
$$
(\id,\iota):(X,\O_X)\lra (Y,\O_Y).
$$
The following should be   well-known:

\begin{lemma}
\mylabel{categorical quotient}
The above morphism of schemes is a   categorical quotient in $(\Sch/k)$.
Moreover, the formation of the quotient is compatible with flat base-change in the scheme $(Y,\O_Y)$.
\end{lemma}

\proof
First note that the inclusion $\O_Y\subset\O_X$ is invariant with respect to multiplication of scalars $\lambda\in k$,  so the morphism belongs to the
category $(\Sch/k)$. Furthermore, the formation of   kernels and finite intersections for maps between quasicoherent sheaves on schemes
is compatible  with flat base-change, and in particular the formation of $(Y,\O_Y)$ is compatible with flat base-change.

We now 
verify the universal property. Let $(T,\O_T)$ be     scheme endowed with the trivial $G$-action,
and $(f,\varphi):(X,\O_X)\ra (T,\O_T)$ be an equivariant  morphism.
Obviously, there is a unique continuous map $g:Y\ra T$ with $f=g\circ\id$. The trivial $G$-action on $(T,\O_T)$ corresponds
to the zero map $\lieg\ra H^0(T,\O_T)$, and equivariance ensures that $f^{-1}(\O_T)\ra\O_X$ factors over the injection $\O_Y\subset\O_X$.
This gives a unique morphism $(g,\psi):(Y,\O_Y)\ra (T,\O_T)$ of ringed spaces that factors  $(f,\varphi)$.
For each point $a\in X$,   the local map $\O_{T,f(a)}\ra \O_{X,a}$ factors over $\O_{Y,a}$, and it follows that $\psi:\O_{T,g(a)}\ra\O_{Y,a}$ is local.
Thus $(g,\psi)$ is a morphism in the category $(\Sch/k)$, which shows the universal property.
\qed

\medskip
We now revert back to the usual notation, and write $Y=X/G$ for the  quotient of the action $\mu:G\times X\ra X$, with quotient map $q:X\ra Y$.
Clearly, this map  is surjective, $Y$ carries the  quotient topology,   and the set-theoretical image of $\mu\times\pr_2:G\times X\ra X\times X$ equals
the fiber product $X\times_YX$.
By construction,  for each open set $U\subset Y$ and each local section $f\in\Gamma(U,q_*(\O_X))$, we have $f\in \Gamma(U,\O_Y)$ if and only if 
$f\circ\mu=f\circ\pr_2$ as morphisms $G\times q^{-1}(U)\ra\AA^1$. 
Summing up, our categorical quotient is also a \emph{uniform geometric quotient},
in the sense of \cite{Mumford; Fogarty; Kirwan 1993}, Definition 0.7.
The following observation will be useful:

\begin{proposition}
\mylabel{quotients normal}
Suppose that $X$ is integral, with function field $F=\O_{X,\eta}$. 
Then   $\O_{Y,a}=\O_{X,a}\cap (F^\lieg)$ for each point $a\in X$.
Moreover, the scheme $Y$ is normal provided this holds for $X$.
\end{proposition}

\proof
Set $R=\O_{X,a}$, such that $R^\lieg=\O_{Y,a}$. Choose a basis $D_1,\ldots,D_n\in\lieg$, and consider
the resulting commutative diagram with exact rows
$$
\begin{CD}
0	@>>>	R^\lieg	@>>>	R	@>>>	R^{\oplus n}\\
@.		@VVV		@VVV		@VVV\\
0	@>>>	F^\lieg	@>>>	F	@>>>	F^{\oplus n}\\
\end{CD}
$$
where the horizontal maps on the right are given by $s\mapsto (D_1(s),\ldots,D_n(s))$.
The commutativity of the left square gives $R^\lieg\subset R\cap F^\lieg$, and the injectivity of the
vertical map on the right ensures the reverse inclusion, by a diagram chase.
Now suppose that $R$ is normal, and $f\in F^\lieg$ satisfies an integral equation over the subring $R^\lieg$.
This is also an integral equation over $R$, hence $f\in R\cap F^\lieg=R^\lieg$.
\qed

\section{Inertia and Jacobson correspondence}
\mylabel{Inertia}

The goal of this section is provide a new, more geometric interpretation of the \emph{Jacobson Correspondence}
(\cite{Jacobson 1937} and \cite{Jacobson 1944}).
We start by recalling this correspondence, which   relates certain subfields and restricted Lie algebras,
in Bourbaki's formulation  (\cite{A 4-7}, Chapter V, \S13, No.\ 3, Theorem 3):

Let $F$ be a field of characteristic $p>0$. It comes with a subfield $F^p$ and
a restricted Lie algebra $\lieg=\Der(F)$ over $F^p$ that is also endowed with the structure of an $F$-vector space.
Note that the bracket is $F^p$-linear but in general not $F$-bilinear. Rather, we have the formula
\begin{equation}
\label{bracket formula}
[\lambda D, \lambda'D']= \lambda\lambda' \cdot [D,D'] + \lambda D(\lambda')\cdot D' - \lambda'D'(\lambda)\cdot D.
\end{equation}
Throughout, a subgroup $\lieh\subset \lieg$ is called an \emph{$F^p$-subalgebra with $F$-multiplication} if it is stable under
bracket, $p$-map, and multiplication by scalars $\lambda\in F$. It is thus a restricted Lie algebra over $F^p$, endowed with
the $F$-multiplication as \emph{additional structure}.
Consider the ordered sets
\begin{gather*}
\Phi=\{ E\mid \text{$F^p\subset E\subset F$  is an intermediate field} \},\\
\Psi=\{\lieh\mid \text{$\lieh\subset \lieg$   is an $F^p$-subalgebra with $F$-multiplication}\}.
\end{gather*}
Similar to  classical Galois theory for separable algebraic extensions, one has inclusion-reversing maps
$\Phi\ra \Psi$ and $\Psi\ra \Phi$
given by 
$$
E\longmapsto \Der_E(F)\quadand   \lieh\longmapsto F^\lieh,
$$
respectively. Here $F^\lieh$ denotes the intersection of the kernels for  $D:F\ra F$, where $D\in \lieh$
runs over all elements. Then the Jacobson Correspondence asserts that the above maps induces a bijection
between the intermediate fields $F^p\subset E\subset F$ having $[F:E]<\infty$ and the $F^p$-subalgebras 
with $F$-multiplication $\lieh\subset\lieg$ having $\dim_F(\lieh)<\infty$.
Moreover, under this bijection   $[F:E] = p^{\dim_F(\lieh)}$ holds.

In particular, if $F$ has \emph{finite $p$-degree}, which means that $F^p\subset F$ is finite,
we get an unconditional identification 
$$
\{\text{intermediate fields $E$} \}= \{\text{$F^p$-subalgebras $\lieh$ with $F$-multiplication}\}.
$$
Forgetting the $F$-multiplication, the restricted Lie algebra  $\lieh=\Der_E(F)$  corresponds to  a height-one group scheme $H$,  with $\lieh=\Lie(H)$.
By construction, this coincides with the Frobenius kernel of the affine group scheme $\Aut_{F/E}$.

We now consider the following set-up    geared towards geometric applications: Let $k$ be a ground field of characteristic $p>0$,
and $F$ be some extension field; one should think of the function field  of some proper integral scheme.
Let $H$ be a height-one group scheme over $k$, with corresponding restricted Lie algebra 
$\lieh=\Lie(H)$. Suppose we have a faithful  action of the group scheme  $H$ on the scheme $\Spec(F)$, in other words, a  homomorphism $\lieh\ra\Der_k(F)$
that is $k$-linear and injective. Throughout, we regard this  homomorphisms also as an inclusion.

Let $E=F^\lieh$,  such that  $\lieh\subset\Der_E(F)$.
Then the field $E$ contains the composite    $k\cdot F^p$, and its spectrum
is the categorical  quotient $\Spec(F)/H$, according to Proposition \ref{categorical quotient}. Moreover, we obtain  subspace $\lieh\subset  \lieh\cdot E\subset  \lieh\cdot F$ inside $\Der_E(F)$.
These are subvector spaces over $k$ and  $E$ and $F$, respectively.  Obviously we have
$$
\dim_F( \lieh\cdot F)\leq \dim_E( \lieh\cdot E)\leq \dim_k(\lieh).
$$
Let us unravel how these various fields and vector spaces are related:

\begin{proposition}
\mylabel{jacobson}
In the above situation, the following holds:
\begin{enumerate}
\item 
The subspace $\lieh\subset \Der_E(F)$ contains an $F$-basis, such that $\lieh\cdot F=\Der_E(F)$.
\item
The  canonical inclusions $E=F^\lieh\subset F^{\lieh\cdot E}\subset F^{\lieh\cdot F}$ are equalities.
\item
The subspace $\lieh\cdot E\subset \Der_E(F)$ is stable with respect to bracket and $p$-map.
\item
The extension $E\subset F$ is finite, of degree $[F:E]=p^{\dim_F( \lieh\cdot F)}$.
\end{enumerate}
\end{proposition}

\proof
To see (ii), choose a $k$-generating set $D_1,\ldots,D_n\in \lieh$.
Clearly, $F^\lieh$ coincides with the intersection of the $\Kernel(D_i:F\ra F)$.
Since $D_1,\ldots,D_n\in \lieh\cdot F$ is an $F$-generating set as well, this intersection
coincides with $F^{\lieh\cdot F}$, and the equalities $F^\lieh =  F^{\lieh\cdot E}=F^{\lieh\cdot F}$ follow.

We next verify that the  $F$-vector subspace $\lieh\cdot F\subset\Der_E(F)$ is stable under bracket and $p$-map.
The former follows from  \eqref{bracket formula}. The latter is then a consequence of the Hochschild  Formula (\cite{Hochschild 1955}, Lemma 1)
$$
(vu)^p = v^pu^p + \ad_{vu}^{p-1}(v)u,
$$
which holds for any $u$ from an associative $\FF_p$-algebra $U$ and $v$ from an $\ad_u$-stable commutative subalgebra $V$.
In turn, $\lieh\cdot F\subset \Der_E(F)$ is an $F^p$-subalgebra
with $F$-multiplication, obviously of finite $F$-dimension. Now the Jacobson Correspondence applied to $E=F^{\lieh\cdot F}$ shows (iv). 
Applying the  correspondence once more reveals $\lieh\cdot F=\Der_E(F)$, and (i) follows.
The above reasoning   likewise shows that the $E$-vector subspace $\lieh\cdot E\subset\Der_E(F)$ is stable under bracket and $p$-map,
which reveals (iii).
\qed

\medskip
We now seek a   more  geometric understanding of the above facts.
Set $\lieh_E=\lieh\otimes_kE$, and consider the $E$-linearization $\lieh_E\ra\Der_E(F)$
of our inclusion   $\lieh\subset\Der_E(F)$.
Write $\lieh_E^\trv\subset\lieh_E$  for  the kernel. This is an ideal,
giving an inclusion $\lieh_E/\lieh_E^\trv\subset\Der_E(F)$.
Now recall that $H$ denotes the height-one group scheme  
with $\Lie(H)=\lieh$. Write $H_E=H\otimes_kE$ for its base-change,
and $H_E^\trv\subset H_E$ for the normal subgroup scheme corresponding to $\lieh_E^\trv$.
This acts trivially on $\Spec(F)$, whereas the quotient $H_E/H_E^\trv$ acts faithfully.

\begin{proposition}
\mylabel{transitive}
The action of the group scheme  $H_E$ on  $\Spec(F)$ is transitive.
\end{proposition}

\proof
Recall that for any site $\catC$, the action of a group-valued sheaf $G$ on a sheaf $Z$ is called
\emph{transitive} if the morphism $\mu\times\pr_2:G\times Z\ra Z\times Z$ is an epimorphism,
where $\mu:G\times Z\ra Z$ denotes the action.

In our situation the site is $(\Aff/E)$, endowed with the fppf topology.
Set $G=H_E/H^\triv_E$ and $Z=\Spec(F)$. We have to check that for any $R$-valued points $a,b\in Z(R)$, 
there is an fppf extension $R\subset R'$
and some $\sigma\in G(R')$ that sends the base-change $a\otimes R'$ to $b\otimes R'$.
Replacing $R$ by $R\otimes_EF$, we may assume that $R$ is an $F$-algebra.
Choose a $p$-basis for the extension $E\subset F$, such that
$F=E[T_1,\ldots, T_r]/(T_1^p-\mu_1,\ldots, T_r^p-\mu_r)$
for some scalars $\mu_i\in E$. Then
$$
F\otimes_ER = R[s_1,\ldots,s_r]/(s_1^p,\ldots,s_r^p)
$$
for the elements $s_i=T_i\otimes 1 - 1\otimes T_i$.
The  $R$-valued points of $Z$ thus correspond to $s_i\mapsto \lambda_i$, where $\lambda_i\in R$ satisfy $\lambda_i^p=0$.
It suffices to treat the case that $a,b\in Z(R)$ is given by   $s_i\mapsto 0$ and $s_i\mapsto \lambda_i$, respectively.
 
The differentials $dT_i\in \Omega^1_{F/E}$ form an $F$-basis. The dual basis inside $\Der_E(F) =\Hom(\Omega^1_{F/E},F)$
are the partial derivatives $ \partial/\partial T_i$. Clearly we have $[\partial/\partial T_i,\partial/\partial T_j]=(\partial/\partial T_i)^\pmap=0$.
Consequently, the linear combination $D=\sum\lambda_i\partial/\partial T_i$ satisfies $D^\pmap=0$, thus $D$ 
is an \emph{additive} element inside $\Der_R(F\otimes_ER)$.
Note that this would fail with coefficients from $F\otimes_ER$ rather than $R$.
By the Demazure--Gabriel Correspondence, it yields a homomorphism of group schemes $\alpha_{p,R}\ra \Aut_{F/E}\otimes_ER$.

According to Proposition \ref{jacobson} we have
$\Der_E(F)=\lieh\cdot F$, so there are elements $D_1,\ldots, D_r\in \lieh\cdot E=\lieh_E/\lieh_E^\triv$ that form an $F$-basis of $\Der_E(F)$.
In particular, we may write $\sum\lambda_i\partial /\partial T_i=\sum \alpha_i D_i$ for some $\alpha_i\in R$.
In turn, we get an additive element  $D\in (\lieh_E/\lieh_E^\triv)\otimes_ER$,
so our homomorphism of group schemes has a factorization $\alpha_{p,R}\ra G_R$.
For $R'=R[\sigma]/(\sigma^p)$ we get a canonical element $\sigma\in\alpha_{p,R'}$,
whose image is likewise denoted by $\sigma\in G(R')$.  By construction, we have 
$$
\sigma^*(s_j) = D(s_j) = \sum_i\lambda_i\partial T_j /\partial T_i  = \lambda_j,  
$$
for all $1\leq j\leq n$, and the desired property $\sigma\cdot a=b$ follows.
\qed

\medskip
Note that   the $E$-scheme $Z=\Spec(F)$ does not contain  a rational point, except for $\lieh=0$.
The existence of such a point would allow
us to form the inertia subgroup scheme and view $Z$ as a \emph{homogeneous space}.
However, we can achieve this    after   further base-change:

Regard $A=F\otimes_EF$ as an $F$-algebra via $\lambda\mapsto 1\otimes \lambda$.
Then the multiplication map $\lambda\otimes\mu\mapsto \lambda\mu$ yields a canonical retraction.
Indeed, $A$ is a local Artin ring with residue field $A/\maxid_A=F$.
In turn, $Z_F=Z\otimes F$ has a unique rational point $z_0\in Z_F$.
Write $H_F^\inert=I(z_0)$ for the resulting \emph{inertia subgroup scheme} inside $H_F=H\otimes_kF$.
By the Demazure--Gabriel Correspondence, it is given by a Lie subalgebra $\lieh_F^{\text{inert}}$ inside $\lieh_F=\lieh\otimes_kF$,
which we call the \emph{inertia  Lie algebra}. 
We now interpret the base change $Z_F$  as a  homogeneous space:

\begin{proposition}
\mylabel{orbit map}
The orbit morphism $H_F\cdot\{z_0\}\ra Z_F$ induces an identification
$H_F/H^\inert_F=\Spec(F\otimes_EF)$. Moreover, the inertia Lie algebra $\lieh_F^\inert$
is the kernel for the canonical surjection 
$$
 \lieh\otimes_k F\lra \lieh\cdot F=\Der_E(F).
$$
Finally, the degree of the field extension $E\subset F$ can be expressed as 
$[F:E]=p^c$, where $c\geq 0$ is the codimension  
of the inertia Lie algebra $\lieh_F^\inert\subset\lieh_F$.
\end{proposition}

\proof
According to Proposition \ref{transitive}, the $H_F$-action on $Z_F$ is transitive, and it follows that
the orbit $H_F\cdot\{z_0\}\ra Z_F$ is an epimorphism. By definition of the inertia subgroup scheme,
the induced   $H_F/H^\inert_F\ra Z_F$ is a monomorphism. Hence the latter is an isomorphism.
This is a finite scheme, and the $F$-dimension for the ring of global sections for the homogeneous space
is given by $p^c$. It follows $[F:E]=p^c$.

It remains to see that the inertia Lie algebra $\lieh_F^\inert$ coincides with the kernel $K$ of
the canonical surjection $\lieh_F\ra \lieh\cdot F$.
We saw in Proposition \ref{jacobson} and the preceding paragraph   that
$$
p^{\dim_F(\lieh\cdot F)} = [F:E]=h^0(\O_{Z}\otimes_EF) = p^{\dim(\lieh_F/\lieh_F^\inert)}.
$$
It thus suffices to verify that the canonical map  $\lieh_F^\inert\ra \Der_E(F)$ is zero.
Suppose this is not the case, and fix some non-zero $D\in\lieh_F^\inert$ with non-zero image.
Choose a $p$-basis for $E\subset F$ and write 
$$
F=E[T_1,\ldots, T_r]/(T_1^p-\mu_1,\ldots, T_r^p-\mu_r)
$$
for some scalars $\mu_i\in E^\times$. The partial derivatives $\partial/\partial T_i\in \Der_E(F)$
form another $F$-basis, and $D=\sum\lambda_i\partial/\partial T_i$. Without restriction,
we may assume $\lambda_1\neq 0$.
Now make a base-change to $R=F$, such that 
$$
A=F\otimes_EF = R[s_1,\ldots,s_r]/(s_1^p,\ldots, s_r^p)
$$
as in the proof for Proposition \ref{transitive}.  Then $D(s_1) =  \lambda_1\otimes 1\not\in\maxid_A$.
But this implies that $H_F^\inert$ does not fix the closed point $z_0\in Z_F=\Spec(A)$, contradiction.
\qed

\section{The foliation rank}
\mylabel{Foliation rank}

Throughout this section, $k$ is a ground field of characteristic $p>0$,
and $X$ is a proper scheme. Note that everything carries over verbatim to proper algebraic spaces.
Let $H=\Aut_{X/k}[F]$ be the resulting height-one group scheme,
whose restricted Lie algebra is $\lieh=H^0(X,\Theta_{X/k})$.
To simplify exposition, we  also  assume   that $X$ is integral. Let 
$\eta\in X$ be the generic point and  $F=k(X)$ be the function field. 
The quotient $Y=X/H$ is integral as well,
and we denote its function field by  $E=k(Y)$.
This field extension $E\subset F$ is finite and purely inseparable. This yields
a numerical invariant:

\begin{definition}
\mylabel{foliation rank}
The \emph{foliation rank} of the proper integral scheme $X$ is the integer $r\geq 0$ defined by the
formula $\deg(X/Y)=p^r$.
\end{definition}

In other words, we have $[F:E]=p^r$. Since the field extension $E\subset F$ has height one, the foliation rank $r\geq 0$
is also given by $r=\dim_F(\Omega^1_{F/E})$, which can also be seen as the
rank of the coherent sheaf $\Omega^1_{X/Y}$.
Dualizing the surjection $\Omega^1_{X/k}\ra\Omega^1_{X/Y}$ gives
an inclusion
$\shF=\Theta_{X/Y}\subset \Theta_{X/k}$.
The subsheaf $\shF$ is closed under Lie brackets and $p$-maps, hence constitutes
a \emph{foliation}, where the integer $\rank(\shF)=\rank(\Omega^1_{X/Y})$ coincides with our foliation rank $r\geq 0$.

To obtain an interpretation  of the foliation rank in terms of group schemes, consider the  restricted Lie algebras 
$$
\lieh=\Lie(H)=H^0(X,\Theta_{X/k}) \quadand \lieg=\Der_E(F)=\Theta_{X/Y,\eta}.
$$
The former is   finite-dimensional over the ground field $k$.
The  latter is finite-dimen\-sional over the function field $E$, and can be seen as the Lie algebra
for the automorphism group scheme for  $\Spec(F)$ viewed as a finite $E$-scheme.
The localization map  
$\lieh=H^0(X,\Theta_{X/k})\ra \Theta_{X/k,\eta} $
respects brackets and $p$-powers, and factors over the subalgebra $\lieg=\Theta_{X/Y,\eta}$.
This gives a   $k$-linear map $\lieh\ra\lieg$, 
together with its $E$-linearization
$$
\lieh\otimes_kE\lra\lieg,\quad \delta\otimes \lambda\longmapsto (f\mapsto \lambda\delta_\eta(f)).
$$
The latter is a homomorphism of restricted Lie algebras over $E$. 
The map $\lieh\ra\lieg$ is injective, because the coherent sheaf $\Theta_{X/k}$ is torsion free, and
we often view it as an inclusion $\lieh\subset\lieg$.
Note, however, that its $E$-linearization  in general it is \emph{neither injective
nor surjective}. This is perhaps the main difference to the classical situation of group actions
rather than group scheme actions.

We are now in the situation studied in Section \ref{Inertia}.
Let $\lieh_F^\inert$ be the inertia Lie algebra inside the base-change $\lieh_F=\lieh\otimes_kF$, corresponding
to the inertia group scheme with respect to the $F$-rational point in $\Spec(F\otimes_EF)$.
From Proposition \ref{orbit map} we obtain:

\begin{proposition}
\mylabel{codimension}
The  foliation rank  $r\geq 0$ of  the scheme  $X$ coincides with the codimension of $\lieh_F^\inert\subset\lieh_F$.
\end{proposition}
 
\medskip
In some sense, this measures how free the Frobenius kernel of the automorphism group scheme
acts generically:

\begin{proposition}
\mylabel{boundary cases}
The foliation rank of the scheme $X$ satisfies $0\leq r\leq h^0(\Theta_{X/k})$.
We have $r=0$ if and only if the Frobenius kernel $H=\Aut_{X/k}[F]$ vanishes.
The condition $r=h^0(\Theta_{X/k})$ holds if and only if $H$ acts freely on some dense open set $U\subset X$.
\end{proposition}

\proof
The inequality  $r\leq h^0(\Theta_{X/k})$ follows from Proposition \ref{codimension}.
If the group scheme $H$ is trivial we have $h^0(\Theta_{X/k})=0$ and hence $r=0$.
Conversely, if $H$ is non-trivial there is a non-zero derivation $D:\O_X\ra\O_X$.
Since the structure sheaf is torsion-free, the derivation remains non-zero at the generic point,
which implies that $E=F^\lieh$ does not coincide with $F$, and thus $r>0$.

If $H$ acts freely on some dense open set, the projection $\epsilon:X\ra Y$ to the quotient $Y=X/H$
is a principal homogeneous $H$-space over the dense open set $V\subset Y$ corresponding to $U$.
In turn $[F:E]=h^0(\O_H)$, and thus $r=\dim(\lieh)=h^0(\Theta_{X/k})$.

Finally, suppose that the foliation rank takes the maximal possible value $r=h^0(\Theta_{X/k})$.
Then the inertia Lie algebra $\lieh_F^\inert\subset\lieh_F$ has codimension $r=\dim(\lieh_F)$,
thus is trivial. It follows that  the group scheme $H_E$ acts freely on $\Spec(F)$ viewed as an $E$-scheme.
Thus there is an open dense set $V\subset Y$ over which the projection $\epsilon:X\ra Y$ 
becomes a principal homogeneous $H$-space, and the $H$-action on $U=\epsilon^{-1}(V)$ is free.
\qed

\medskip
We next describe how the foliation rank behaves under   birational maps:

\begin{proposition}
\mylabel{birational}
Let $f:X\ra X'$ be a birational morphism to another proper integral scheme $X'$, with the property $\O_{X'}=f_*(\O_X)$.
Then the respective foliation ranks satisfy $r\leq r'$. 
\end{proposition}

\proof
According to Blanchard's Lemma, there is a unique 
homomorphism $f_*:\Aut_{X/k}^0\ra \Aut_{X'/k}^0$ of group schemes making the morphism $f:X\ra X'$ equivariant.
Indeed, the original form of the lemma for complex-analytic spaces (\cite{Blanchard 1956}, Proposition I.1) was extended  to schemes
by Brion, Samuel and Uma (\cite{Brion; Samuel; Uma 2013}, Proposition 4.2.1).

The homomorphism of group schemes is a monomorphism, because $f$ is birational, and the schemes in question are integral.
In particular, the induced homomorphism on Frobenius kernel gives a closed embedding $H\subset H'$,
and an injection $\lieh\subset \lieh'$ of restricted Lie algebras.
For the common function field $F=k(X)=k(X')$, we get $F^\lieh\supset F^{\lieh'}$, and $r\leq r'$ follows.
\qed

\medskip
The following gives an upper bound on the foliation rank:

\begin{proposition}
\mylabel{upper bound}
Let $i\geq 0$ be some integer, and suppose that 
 the   coherent sheaf $\shF=\uHom(\Omega^i_{X/k},\O_X)$ satisfies $h^0(\shF)=0$. Then $X$ has foliation rank $r<i$.
\end{proposition}

\proof
We have to show that the vector space  $\lieh\cdot F=\Der_E(F)$ has dimension at most $i-1$.
Seeking a contradiction, we suppose that there
are   $k$-derivations $D_1,\ldots,D_i:\O_X\ra \O_X$ that are $F$-linearly independent.
Then the same holds for the corresponding $\O_X$-linear maps $s_1,\ldots,s_i:\O_X\ra\Theta_{X/k}$.
Consequently their wedge product $s_1\wedge\ldots\wedge s_i:\O_X\ra \Lambda^i(\Theta_{X/k})$ is generically non-zero.
The universal property of exteriour powers gives a canonical map
$\Lambda^i(\Theta_{X/k})\ra\uHom(\Omega^i_{X/k},\O_X)=\shF$, which  is generically bijective.
Thus $s_1\wedge\ldots\wedge s_i$ yield a non-zero global section of $\shF$, contradiction.
\qed

\medskip
Recall that our proper integral   $X$ comes with a \emph{dualizing sheaf} $\omega_X$ and  a \emph{trace map} $H^n(X,\omega_X)\ra k$,  such that
the  ensuing  pairing $\Hom(\shF,\omega_X)\times H^n(X,\shF)\ra k$ is non-degenarate. Here $n=\dim(X)$ and $\shF$ is coherent.  

\begin{corollary}
\mylabel{normal surface}
Let $X$ be a geometrically normal surface  with $h^0(\omega^{\vee}_X)=0$. 
Then the foliation rank is $r\leq 1$.
\end{corollary}

\proof
Replacing the ground field $k$ by the field $H^0(X,\O_X)$, it suffices to treat the case $h^0(\O_X)=1$.
By Serre's Criterion, the scheme $X$ is regular in codimension one, so the locus of non-smoothness $\Sing(X/k)$ is finite.
Let $f:S\ra X$ be a resolution of singularities.
Suppose for the moment that the regular surface $S$ is smooth. Then $\omega_S=\Omega^2_{S/k}$.
Consider the following chain of canonical maps
$$
\Omega^2_{X/k}\lra f_*f^*(\Omega_{X/k}^2)\lra f_*(\Omega^2_{S/k}) \lra f_*(\omega_S)\lra \omega_X,
$$
where to the right is the trace map.
All these maps are bijective on the complement $U=X\smallsetminus \Sing(X/k)$, so the same holds for the dual map 
$$
\varphi:\omega_X^\vee\lra \uHom(\Omega^2_{X/k},\O_X)=\shF.
$$
According to \cite{Hartshorne 1994}, Corollary 1.8 and Theorem 1.9, these rank-one sheaves are reflexive
and satisfy the Serre Condition $(S_2)$.
Since $\varphi|U$ is bijective, already $\varphi$ is bijective, by loc.\ cit.\ Theorem 1.12.
The assertion thus follows from the theorem.

It remains to treat the case that the ground field $k$ is imperfect. Choose a perfect closure $k'$.
The base-change $X'=X\otimes_kk'$ is normal, and the above reasoning applies to any resolution of singularities
$S'\ra X'$. It follows that $\omega_X^\vee$ and $\shF$ become isomorphic after base-changing to $k'$.
If follows that $\Hom(\omega_X,\shF)$ is one-dimensional.
Choose a non-zero element $\varphi:\omega_X\ra \shF$. Then $\varphi\otimes k'$ must be bijective,
and by descent the same holds for $\varphi$.
\qed

\medskip
This applies in particular to smooth surfaces $S$ of Kodaira dimension $\kod(S)\geq 1$,
which comprise surfaces of general type, and the properly elliptic surfaces, including those with  quasi-elliptic fibration.
It also applies to surfaces $S$ with Kodaira dimension zero, provided that
the dualizing sheaf of the minimal model $X$ is non-trivial.

Let $S$ be a smooth surface of general type, and $X$ be its canonical  model.
This is the homogeneous spectrum $P(S,\omega_S)$ of the graded ring $R(S,\omega_S)=\bigoplus H^0(S,\omega_S^{\otimes t})$.
Then $X$ is normal, the singularities are at most rational double points, and the dualizing sheaf $\omega_X$ is ample.
We also say that $X$ is a \emph{canonically polarized surfaces}.
Obviously   $h^0(\omega_X^{\otimes-1})=0$, and   $X$ has foliation rank $r\leq 1$.
According to Proposition \ref{birational}, the same holds for $S$.

\begin{proposition}
\mylabel{local description}
Suppose that $X$ has foliation rank $r=1$, and let $D\in H^0(X,\Theta_{X/k})$ be any non-zero global section.
Then for each point $x\in X$, the local ring $\O_{Y,\epsilon(x)}$ is the kernel for the additive map
$D:\O_{X,x}\ra\O_{X,x}$.
\end{proposition}

\proof
Set $y=\epsilon(x)$. The local ring is given by $\O_{Y,y}=\O_{X,x}^\lieh$, which is contained in the kernel
$\O_{X,x}^D$ of the derivation $D$. Let $f\in \O_{X,x}^D$, and $D'\in\lieh$ be another derivation.
Then $D'=\lambda D$ for some element $\lambda$ from the  function field $F=\Frac(\O_{X,x})$, and thus $D'(f) =\lambda D(f) =0$ inside $F$.
Since the localization map $\O_{X,x}\ra F$ is injective,
we already have $D'(f)=0$ inside $\O_{X,x}$.
This shows $f\in \O_{Y,y}$. In turn, the inclusion $\O_{Y,y}\subset\O_{X,x}^D$ is an equality.
\qed

\medskip
We will later see that for $r=1$ each vector in $\lieg$ is $p$-closed. Thus the non-zero elements
$D\in \lieh$ indeed yield   height-one group schemes $N\subset H$ of order $|N|=p$,
such that  $Y=X/N$.

\section{Invariant   subspaces}
\mylabel{Invariant subspaces}

Let $k$ be  a ground field of characteristic $p\geq 0$ and $V$ be a finite-dimensional vector space of dimension $n\geq 0$.
Let us write $\GL_{V/k}$ for the group-valued functor on the category $(\Aff/k)$ of affine $k$-schemes $T=\Spec(R)$
defined by
$$
\GL_{V/k}(R) = \Aut_R(V\otimes_kR).
$$
This satisfies the sheaf axiom with respect to the fppf topology.
In fact, it is representable by an affine group scheme, and the choice of a basis $e_1,\ldots,e_n\in V$ yields  $\GL_{V/k}\simeq\GL_{n,k}$.

Let us  write $\underline{V}$ for the abelian functor whose group of  $R$-valued points  is 
$\underline{V}(R)=V\otimes_kR$. As explained in \cite{EGA I}, Chapter I, Section 9.6 this is  represented by an affine scheme, namely 
the spectrum of the symmetric algebra on the dual vector space $V^*$. Moreover, the structure morphism $\underline{V}\ra\Spec(k)$ carries the
structure of a \emph{vector bundle} of rank $n$ with $\underline{V}(k)=V$, and the canonical homomorphism
$
\GL_{V/k}\ra \Aut_{\underline{V}/k}
$
of group schemes is bijective. Combining \cite{GB III}, Theorem 11.7 with \cite{SGA 4b}, Expos\'e VIII, Corollary 2.3 and \cite{A 4-7}, Chapter V, \S10, No.\ 5, Proposition 9
one sees that each $\GL_{V/k}$-torsor  is trivial,
that is, the non-abelian cohomology set $H^1(k,\GL_{V/k})$ with respect to the fppf topology is a singleton. In other words,   all vector bundles $E\ra\Spec(k)$
of rank $n$ are isomorphic to $\underline{V}$. 

Now let $H\subset\GL_{V/k}$ be a   subgroup scheme, and    $T\ra\Spec(k)$ be a $H$-torsor. Then the quotient
$$
{}^T\underline{V} = H\backslash (T\times\underline{V}) = T\wedge^H\underline{V}
$$
with respect to the diagonal action $\sigma\cdot (t,v)= (\sigma t,\sigma v)$
is another vector bundle called the \emph{$T$-twist}. Note that under the identification of left and right action,
the above  action can also be viewed as $\sigma\cdot (t,v)=(t\sigma^{-1},\sigma v)$, which explains the notation
$T\wedge^H\underline{V}$.
\emph{We now consider the following general problem:
What subbundles exist in the $T$-twist whose pull-back to $T$ are contained in the pullback
of a fixed subbundle $\underline{V}'\subset \underline{V}$?} By fppf descent, these pullbacks
correspond to subbundles inside the induced bundle $\underline{V}\times T\ra T$ whose total space 
is invariant with respect to the diagonal $H$-action.

The $n$-dimensional vector space  ${}^TV=({}^T\underline{V})(k)$  of $k$-rational points is likewise called
the \emph{$T$-twist} of $V$. If $H$ is finite and $T=\Spec(L)$ is the spectrum of a field, we are thus looking
for $k$-vector subspaces $U\subset {}^TV$ such that  $U\otimes_kL$ is contained in the base-change
$V'\otimes_kL$, or equivalently to $L$-vector subspaces
in $V'\otimes_kL$ that are invariant for the diagonal $H$-action.

Suppose now that $p>0$, and that $H=\alpha_p$ is the \emph{infinitesimal group scheme}
defined by $H(R)= \{\alpha\in R\mid \alpha^p=0\}$, where the group law is given by addition.
Recall that the Lie algebra of $\GL_{V/k}$ is the vector space $\liegl(V)=\End_k(V)$,
where the Lie bracket  is given by commutators $[f,g]=fg-gf$, and the $p$-map 
$f^\pmap=f^p$ is the $p$-fold composition. 
The inclusion homomorphism $H\ra \GL_{V/k}$ corresponds to a vector 
$f\in\liegl(V)$ that is   nilpotent, with
all   Jordan blocks of size $\leq p$. On $R$-valued points, the map becomes
$$
H(R)\lra \GL_{V/k}(R),\quad \alpha\longmapsto  \sum_{i=0}^{p-1} \frac{(\alpha f)^i}{i!}.
$$
Set $ e^{\alpha f} =  \sum_{i=0}^{p-1}(\alpha f)^i/i!$ to simplify notation.
By naturality, the above maps are determined by the single matrix $e^{t f}$
with entries in the truncated polynomial ring  $R=k[t]/(t^p)$.
The following   is well-known:

\begin{lemma}
\mylabel{torsor description}
Each  torsor $T$ for the infinitesimal group scheme $H=\alpha_p$ is isomorphic to the spectrum of     $L=k[s]/(s^p-\omega)$ for some $\omega\in k$,
where the group elements $\alpha\in H(R)$ act via $s\mapsto s+\alpha$.
The torsor $T$ is non-trivial if and only if $L$ is a field. Moreover, for
each purely inseparable field extension $k\subset L$ of degree $p$, the spectrum $\Spec(L)$ admits the structure of 
a $H$-torsor.
\end{lemma}

\proof
Consider the relative Frobenius map $F:\GG_a\ra\GG_a$ on the additive group, which comes
from the $k$-linear map $k[t]\ra k[t]$ given by $t\mapsto t^p$. Then $H=\alpha_p$ is the kernel.
The short exact sequence $0\ra H\ra\GG_a\stackrel{F}{\ra} \GG_a\ra 0$ yields a long exact sequence
$$
k\lra k \lra H^1(k,H)\lra H^1(k,\GG_a)\lra H^1(k,\GG_a).
$$
The terms on the right vanish. It follows that each $H$-torsor $T$ arises as the fiber
for $F:\GG_a\ra\GG_a$ over some rational point $\omega\in\GG_a(k)$.
Thus $T$ is equivariantly isomorphic to the  spectrum of  $k[s]/(s^p-\omega)$, 
where    the group elements $\alpha\in H(R)$ act via $s\mapsto s+\alpha$.
If $T$ is non-trivial, the polynomial $s^p-\omega\in k[s]$ has no root in $k$.
We infer that it is irreducible, because  the algebra $L=k[s]/(s^p-\omega)$ has prime degree $p$.
Thus $L$ is a field, which is purely inseparable over $k$. Conversely, if $L$ is a field,
then $T$ has no rational point, and the torsor is non-trivial.

Finally, let $k\subset L$ be a purely inseparable extension of degree $p$. For each element in $L$ not contained in $k$ we get
an identification $L=k[s]/(s^p-\omega)$. Thus $\Spec(L)$ arises as fiber of the relative Frobenius map, thus admits the structure
of a $H$-torsor. 
\qed

\medskip
For the applications we have in mind, we now consider the particular situation  that $V=k[t]/(t^p)$ is 
the underlying vector space of dimension $n=p$ coming   from the 
truncated polynomial ring, and $V'=tk[t]/(t^p)$ is given by the maximal ideal. 
Each vector can be uniquely written as a polynomial $f(t)=\sum_{i=0}^{p-1} \lambda_i t^i$,
with coefficients $\lambda_i\in k$. This vector space comes with a canonical action of the additive group $\GG_a$,
where the elements $\alpha\in \GG_a(R)=R$  act  via  $f(t)\mapsto f(t+\alpha)$.
With respect to the canonical basis $t^0,\ldots,t^{p-1}\in V\otimes_kR$, this automorphism is given by $\alpha\mapsto (\alpha_{ij})$,
where the matrix entries are $\alpha_{ij}=\binom{j}{j-i}\alpha^{j-i}$. In turn, we get an induced action of the Frobenius kernel
$H=\alpha_p$.
Note that $V'\subset V$ is not $H$-invariant, because some $\alpha_{0j}=\alpha^j$ are non-zero for $\alpha\neq 0$.

Now let $T=\Spec(L)$ be a   $H$-torsor.
The resulting twist ${}^TV$  is another vector space of dimension $n=p$.
Note that both $V$ and ${}^TV$ are isomorphic to $k^{\oplus p}$, but there is no canonical isomorphism.
The following observation will be crucial for later applications:

\begin{proposition}
\mylabel{no subspace}
In the above situation, there is no vector $x\neq 0$ inside the twist  ${}^TV$
such that the induced element $x\otimes 1$ inside ${}^TV\otimes_kL=V\otimes_kL$ is contained in the base change $V'\otimes_kL$.
\end{proposition}

\proof
Seeking a contradiction, we assume that such an element exists. Its image $x\otimes 1$ inside ${}^TV\otimes_kL=V\otimes_kL$ takes the form
$f(t)=\sum_{i=1}^{p-1} \lambda_it^i$   with coefficients $\lambda_i\in L$.
According to Lemma \ref{torsor description}, we have  $L=k[s]/(s^p-\omega)$ for some $\omega\in k$,
and the group elements $\alpha\in H(R)$ act via $s\mapsto s+\alpha$. 
Write $\lambda_i=\sum_{j=0}^{p-1}\lambda_{ij}s^j$ with   coefficients $\lambda_{ij}\in k$.
The  $H$-invariance of the vector $f(t)\in V\otimes_kL$ with respect to the diagonal $H$-action means 
\begin{equation}
\label{invariance}
\sum_{i=1}^{p-1}\sum_{j=0}^{p-1}\lambda_{ij}(s+\alpha)^j(t+\alpha)^i= \sum_{i=1}^{p-1}\sum_{j=0}^{p-1}\lambda_{ij}s^jt^i
\end{equation}
for each $\alpha\in H(R)$. Our task is to infer $\lambda_{ij}=0$.
We now consider the universal situation, where $\alpha$ is the class of the
indeterminate in the truncated polynomial ring $R=k[u]/(u^p)$. Then \eqref{invariance} becomes an equation in the residue class  ring
$k[t,s,u]/(t^p,s^p-\omega, u^p)$. Writing
$$
(s+\alpha)^j(t+\alpha)^i= s^jt^i + \alpha(js^{j-1}t^i+is^jt^{i-1})  + \alpha^2(\ldots)
$$
as a polynomial in $\alpha$ and comparing coefficients in \eqref{invariance} at the linear terms   we get the equation
\begin{equation}
\label{double sum}
\sum_{i=1}^{p-1}\sum_{j=0}^{p-1}\lambda_{ij}(js^{j-1}t^i+is^jt^{i-1})=0.
\end{equation}
The following argument, more elegant than  our original reasoning, was indicated by the referee: 
To see that $\lambda_{ij}$ vanishes, it  suffices to check that the polynomial
$$
F=\sum_{i=1}^{p-1}\sum_{j=0}^{p-1}\lambda_{ij} s^jt^i .
$$
is divisible by $t^h$ for   $1\leq h\leq p$ inside the  factorial ring $k[s,t]$.  This is obvious for $h=1$. Suppose now $2\leq h\leq p$,
and that  $F=t^{h-1}G$ for some polynomial $G$. Then
$$
\partial F/\partial s = t^{h-1}\partial G/\partial s\quadand 
\partial F/\partial t =  t^{h-1}\partial G/\partial t + (h-1)t^{h-2}G.
$$
Equation \eqref{double sum} means that $\partial F/\partial s+ \partial F/\partial t=0$. Together with the above computation,
this  gives $t\mid G$ and hence $t^h\mid F$.
\qed

\section{Automorphisms for prime-degree radical extensions}
\mylabel{Radical}

Let $k$ be a ground field of characteristic $p>0$. For each scalar $\omega\in k$, write 
$$
L=L_\omega=k[t]/(t^p-\omega)
$$ 
for the resulting finite algebra of rank $p$. Each element can be uniquely written
as $\sum_{i=0}^{p-1}\lambda_it^i$, and we call such expressions \emph{truncated polynomials}.
Write $\Aut_{L/k}$ for the group-valued functor on the category $(\Aff/k)$ whose
$R$-valued points are the $R$-linear automorphisms of $L\otimes R$.
This functor is representable by an affine group scheme.
In this section we make  a detailed study of  the opposite group  scheme
$$
G=G_\omega=\Aut_{\Spec(L_\omega)/k}=(\Aut_{L_\omega/k})^\op,
$$
which comprises the automorphisms of the affine scheme $\Spec(L)$.
We shall see that $G$ is  non-smooth, so understanding the scheme structure is of paramount importance.
Note that the  Lie algebras $\lieg=\Lie(G)$  were discovered by Witt,
compare the discussions in \cite{Chang 1941}, Introduction  and also \cite{Zassenhaus 1939}, footnote on page 3.
These so-called \emph{Witt algebras} will be studied in the next section.

Any     automorphism $g:L\otimes_k R\ra L\otimes_k R$ is determined by the image of the
generator $t$, which is some truncated polynomial $\varphi_g(t)=\sum_{i=0}^{p-1} \alpha_it^i$.
The multiplication $gh\in G(R)$ of group elements corresponds to the substitution $\varphi_h(\varphi_g(t))$ of truncated polynomials.

The inverse group element $g^{-1}$ defines another truncated polynomial $\varphi_{g^{-1}}=\sum_{i=0}^{p-1}\beta_it^i$,
such that $\sum\alpha_i(\sum\beta_jt^j)^i=t=\sum\beta_i(\sum\alpha_jt^j)^i$. 
Note, however, that the truncated polynomials attached to group elements are never
units in the polynomial ring $R[t]$, unless $R=0$. To avoid this   ambiguity in notation
we use the additional symbol $\varphi_g(t)$ to denote the image of the indeterminate under $g\in G(R)$.

The coefficients in the truncated polynomials $\varphi_g(t)= \sum_{i=0}^{p-1}\lambda_it^i$ for the group elements $g\in G(R)$
define a monomorphism $G \ra\AA^p$.

\begin{proposition}
\mylabel{embedding}
The monomorphism $G \ra \AA^p$ is   an embedding, and its image is the  intersection of the closed set defined
by the Fermat equation
\begin{equation}
\label{equation group}
\lambda_0^p + (\lambda_1-1)^p\omega + \lambda_2^p\omega^2 + \ldots + \lambda_{p-1}^p\omega^{p-1}=0
\end{equation}
with the open set given by $\det(\alpha_{ij})\neq 0$.
Here the matrix entries come from the truncated polynomials $(\sum\lambda_it^i)^j=\sum \alpha_{ij}t^i$, with $0\leq i,j\leq p-1$.
\end{proposition}

\proof
For each $g\in G(R)$, with truncated polynomial $\varphi_g(t)= \sum \lambda_it^i$, the images $(\sum\lambda_it^i)^j$
of the basis vectors $t^j$ form a $R$-basis of $L\otimes R$, thus  $G\ra\AA^p$ factors over
the open set $U\subset\AA^p$ given by $\det(\alpha_{ij})\neq 0$.
Since $t^p=\omega$, we also have $(\sum \lambda_it^i)^p=\omega$, so the monomorphism also factors over the closed set $Z\subset\AA^p$
defined by \eqref{equation group}. Any tuple $(\lambda_0,\ldots,\lambda_{p-1})\in \AA^p(R)$ lying in $U\cap Z$
gives via the truncated polynomial $\sum\lambda_it^i$ some group element $g\in G(R)$.
It follows that the monomorphism $G\ra \AA^p$ is an embedding, with image $U\cap Z$.
\qed

\medskip
Note that  throughout, we  regard the coefficients $\lambda_i$ either as scalars or as indeterminates, depending on the context.
This abuse of notation simplifies   exposition and should not cause   confusion.

\begin{proposition}
\mylabel{neutral element}
The neutral element $e\in G$ has coordinates $(0,1,0,\ldots,0)$ with respect to the embedding $G\subset \AA^p$.
If the scalar $\omega\in k$ is not a $p$-power, then the group of rational points is $G(k)=\{e\}$.
\end{proposition}

\proof
The truncated polynomial of the neutral element is $\varphi_e(t)=t$, which gives the coordinates of $e\in G$.
Now suppose that $\omega\not\in k^p$. By Proposition \ref{embedding}, it suffices to verify that  the polynomial equation 
$\omega^0T_0^p+\omega^1T_1^p+\ldots+\omega^{p-1}T_{p-1}^p=0$ has no non-trivial solution.
The latter means that $1,\omega,\ldots,\omega^{p-1}\in k$ are linearly independent over $k^p$.
This indeed holds, because $k^p\subset k$ is an extension of height $\leq 1$, hence the minimal polynomial
of any $\lambda\in k$ not contained in $k^p$ is of the form $T^p-\lambda^p$.
\qed

\medskip
We now consider the Frobenius pullback $G^\frob$ and its reduced part $G^\frob_\red=(G^\frob)_\red$.
Note that over imperfect fields, reduced parts of group schemes may fail to be subgroup schemes,
see \cite{Fanelli; Schroeer 2020}, Proposition 1.6 for an example.
The following shows that even if it is a subgroup scheme, it might be non-normal. Note that 
this phenomenon seems to be the crucial ingredient for the main results of this paper.

\begin{proposition}
\mylabel{reduced part}
The reduced part  $G^\frob_\red\subset G^\frob$ is a non-normal subgroup scheme. Moreover,  $G^\frob_\red\simeq
U\rtimes\GG_m$, where $U$ has a composition series of length $p-2$  whose quotients are isomorphic to the additive group $\GG_a$.
In particular, $G$ is affine, irreducible, and of dimension $p-1$.
\end{proposition}

\proof
Recall that $G^{(p)}$ is defined as the base-change of $G$ with respect to the absolute Frobenius map on $\Spec(k)$.
Clearly, the Frobenius pullback of $L_\omega=k[t]/(t^p-\omega)$ is isomorphic to $L_0=k[t]/(t^p)$,
and for our automorphism group schemes this means $G^\frob\simeq G_0$.
Thus we may assume $\omega=0$, and work with $G= G^\frob$.

The embedding $G\subset\AA^p$ in Proposition \ref{embedding} is now  given by the conditions $\lambda_0^p=0$ and $\det(\alpha_{ij})\neq 0$.
View the entries of the matrix $(\alpha_{ij})$ as elements from the ring $A=k[\lambda_0,\ldots,\lambda_{p-1}]/(\lambda_0^p)$.
Taken modulo the radical $\Rad(A)=(\lambda_0)$, the matrix takes lower triangular form,
with diagonal entries $1,\lambda_1,\ldots,\lambda_1^{p-1}$. In turn, the embedding $G\subset\AA^p$
is given by $\lambda_0^p=0$ and $\lambda_1\neq 0$.
Consequently, the   reduced part $G_\red$ is defined by $\lambda_0=0$ and $\lambda_1\neq 0$, which is smooth.
Moreover, we see that $G$ is affine, irreducible, and of dimension $p-1$.

Given two truncated polynomials $\varphi_g=\sum\alpha_it^i$ and $\varphi_h=\sum\beta_it^i$ with constant terms $\alpha_0=\beta_0=0$,
the substitution $\varphi_g(\varphi_h(t))$ also has constant term zero,
so the subsets $G_\red(R)\subset G(R)$ are subgroups. Over  $R=k[u,v,\epsilon]/(uv-1,\epsilon^2)$,
the truncated polynomials $\varphi_g=\epsilon + t$ and $\varphi_h=ut$ yield
$$
\varphi_{g^{-1}}(t)=-\epsilon +t \quadand \varphi_{g^{-1}}(\varphi_h(\varphi_g(t)))= \epsilon (u-1)+ ut,
$$
so the subgroup $G_\red(R)\subset G(R)$ fails to be normal. 

Summing up, $G_\red\subset G$ is a smooth non-normal subgroup scheme.
The map $\sum_{i=1}^{p-1}\lambda_it^i\mapsto \lambda_1$ defines a short exact sequence
$$
0\lra U\lra G_\red\lra \GG_m\lra 0.
$$
The inclusion $U\subset\AA^p$ is given by $\lambda_0=0$ and $\lambda_1=1$, hence the underlying scheme of $U$ is a copy
of the affine space $\AA^{p-2}$. By Lazard's Theorem  (\cite{Demazure; Gabriel 1970}, Chapter IV,\S4, 4.1), 
the group scheme $U$ admits a composition series   whose quotients
are isomorphic to the additive group $\GG_a$. Moreover, the projection $G_\red\ra\GG_m$ has
a section  via $\lambda_1\mapsto \lambda_1t$. This is a homomorphism, hence  $G_\red$ is a semidirect product.
\qed

\medskip
Clearly, the group elements $g\in G(R)$ with linear truncated polynomial $\varphi_g=\lambda_0+\lambda_1t$ form a closed subgroup scheme $B\subset G$.
It sits in a short exact sequence
\begin{equation}
\label{borel extension}
0\lra \alpha_p\lra B\lra\GG_m\lra 0,
\end{equation}
where the map on the left is given by $\lambda_0\mapsto \lambda_0+t$, and the map on the right comes from $\lambda_0+\lambda_1t\ra\lambda_1$.
In particular, $B$ is a connected solvable group scheme. We see later that $B$ is maximal with respect to this property,
so one may regard it as a  Borel group. However, we want to stress that this lies in a non-smooth group scheme,
and $B$ itself is non-reduced. We therefore call $B\subset G$ a \emph{non-reduced Borel group}.

An element $\lambda_1\in\GG_m(R)$ lies in the image if and only if $\lambda_0=(1-\lambda_1)\omega^{1/p}$ exists in $R$.
It follows that the extension \eqref{borel extension} splits if and only if $\omega\in k$ is a $p$-power. Moreover, we see that the canonical map
$\GG_m\ra\Aut_{\alpha_p/k}=\GG_m$ is the identity. Note that the group of all such extension of $\GG_m$ by $\alpha_p$, with non-trivial $\GG_m$-action,
is identified with $k/k^p$, according to \cite{Demazure; Gabriel 1970}, Chapter III, \S6, Corollary 6.4.
Note that for $p=2$ the inclusion $B\subset G$ is an equality.
In any case, the pullback of the extension \eqref{borel extension} along the inclusion 
$\mu_p\subset\GG_m$ admits a splitting given by $\lambda_1\mapsto \lambda_1t$,
and one sees $B\times_{\GG_m}\mu_p=\alpha_p\rtimes\mu_p$.

Write $G[F]$ for the kernel of the relative Frobenius map $G\ra G^{(p)}$, which is a normal subgroup scheme of height one.
We now consider the resulting  $G/G[F]\subset G^{(p)}$.

\begin{proposition}
\mylabel{frobenius quotient}
The group scheme $G/G[F]$ is smooth, and coincides with the reduced part $G^\frob_\red$ inside  the Frobenius pullback $G^\frob$.  
\end{proposition}

\proof
We may assume that $k$ is algebraically closed. We first verify that   $G/G[F]$ is reduced.
The short exact sequence \eqref{borel extension} yields an inclusion $\alpha_p\subset G$. This is not normal, 
but contained in the Frobenius kernel $G[F]$. The resulting  projection 
$G/\alpha_p\ra G/G[F]$ is faithfully flat,
and it suffices to check that the homogeneous space  $G/\alpha_p$ is reduced.
Since $G$ acts transitively, it is enough to verify that the 
the local ring at the image in $G/\alpha_p$ of the origin $e\in G$ is regular.
According to \cite{Schroeer 2007}, Proposition 2.2 it is enough to check that   in the local ring 
$\O_{G,e}$, the ideal $\ideala$ corresponding to the subgroup scheme $\alpha_p\subset G$ 
has finite projective dimension.
But this is clear, because it is given by the complete intersection $\lambda_1=\ldots=\lambda_{p-1}=0$.

Thus $G/G[F]$ is reduced. The reduced closed subschemes $G/G[F]$ and $G^\frob_\red$ inside the Frobenius pullback
have the same underlying  set, whence $G/G[F]=G^\frob_\red$. 
The latter is smooth by Proposition \ref{reduced part}, whence the same holds for the former.
\qed

\medskip
Now consider the conjugacy map $c:G\ra\Aut_{G/k}$, sending $g\in G(R)$ to the
automorphism $x\mapsto gxg^{-1}$.  In terms of truncated polynomials,
$gxg^{-1}$ is given by the triple substitution $\varphi_{g^{-1}}(\varphi_x(\varphi_g(t)))$.

\begin{proposition}
\mylabel{conjugacy map}
The conjugacy map $c:G\ra\Aut_{G/k}$ is an isomorphism.
\end{proposition}

\proof
For $\omega=0$, this holds by \cite{Sancho de Salas 2000}, Theorem 4.13. The general case follows by base-changing to $k^\alg$,
and using descent.
\qed

\medskip
In other words, the center    and the
scheme of outer automorphisms are trivial. One also says the the group scheme $G$ is \emph{complete}. 
Now recall that  $G=G_\omega$ depends on a scalar $\omega\in k$. 
 
\begin{proposition}
\mylabel{twisted form}
For each pair of scalars   $\omega,\omega'\in k^\times$, the following are equivalent:
\begin{enumerate}
\item The  $k$-algebras $L_\omega$ and $L_{\omega'}$ are isomorphic.
\item The group schemes $G_{\omega}$ and $G_{\omega'}$   are isomorphic.
\item We have $k^p(\omega)=k^p(\omega')$ as subfields inside $k$.
\end{enumerate} 
\end{proposition}

\proof
According to \cite{Giraud 1971}, Chapter III, Corollary 2.5.2 the category of twisted forms for $L_\omega$,
and the category of twisted forms of $G_\omega$  are both equivalent to the category of $G_\omega$-torsors.
This implies the equivalence of (i) and (ii).

It remains to check (i)$\Leftrightarrow$(iii). Write $L_\omega=k[t]/(t^p-\omega)$
and $L_{\omega'}=k[t']/(t'^p-\omega')$.
Suppose first that these algebras are isomorphic. Choose an isomorphism and regard it
as an identification $L_\omega=L_{\omega'}$. Then $t'=\sum_{i=0}^{p-1}\lambda_it^i$,
consequently $\omega'=\sum\lambda_i^p\omega^i$, and hence $k^p(\omega')\subset k^p(\omega)$.
By symmetry, the reverse implication holds as well.

Conversely, suppose that $k^p(\omega)=k^p(\omega')$. If this subfield coincides with $k^p$,
then both $\omega,\omega'\in k$ are $p$-powers, hence both algebras $L_\omega$, $L_{\omega'}$ are isomorphic
to $k[t]/(t^p)$.  Suppose now that the subfield is different from $k^p$.
Taking $p$-th roots  we get $k(\omega^{1/p})=k(\omega'^{1/p})$ inside some perfect closure $k^\perf$.
These fields are isomorphic to $L_\omega$ and $L_{\omega'}$, because both scalars 
$\omega,\omega'\in k$ are not $p$-powers.
\qed

\medskip
In particular,
each $L=L_\omega$ is a twisted form of $L_0$, and 
each $G=G_\omega$ is a twisted form of $G_0$. Up to isomorphism, these twisted forms correspond to classes 
in non-abelian cohomology set $H^1(k,G_0)$. We will use this throughout to gain insight into $G$,
by using facts on $G_0$.  
For example, from Proposition \ref{embedding} we see that the \emph{locus of non-smoothness} $\Sing(G_0/k)$, defined as in 
\cite{Fanelli; Schroeer 2020}, Section 2, equals the whole scheme $G_0$. Hence the same holds for $G$, because it
is a twisted form of $G_0$.
 
We now write $\Sing(G)$ for the \emph{singular locus} of $G$, which comprise all points   $a\in G$
where the local ring $\O_{G,a}$ is singular. Note that the formation of such loci commutes with 
base-changes along separable extension, but usually not with   inseparable extensions.

\begin{proposition}
\mylabel{singular locus}
The local ring at the origin  is singular, with embedding dimension $\edim(\O_{G,e})=p$.
Moreover, the inclusion $\Sing(G)\subset G$   is not an equality if and only
if $\omega\in k$ is not a $p$-power.  In this case,  the singular locus 
has codimension one in $G$.
\end{proposition}
 
\proof
Since $e\in G$ is a rational point, the embedding dimension of $\O_{G,a}$ does not change under ground field
extensions.  
If $\omega\in k^p$ we have $G\simeq  G_0$ and thus for every point $a\in G$ the local ring $\O_{G,a}$ is singular.
Now suppose that $\omega$ is not a $p$-power, and consider the $p$-Fermat hypersurface
$X\subset\PP^{p-1}$ defined by the homogeneous polynomial $\lambda_0T_0^p+\ldots+\lambda_{p-1}T_{p-1}^p$,
with coefficient $\lambda_i=\omega^i$.  
The field extension $k^p\subset E$ generated by $\lambda_i/\lambda_0=\omega^i$ is nothing but $k^p(\omega)$.
It has degree $[E:k]=p$, hence its $p$-degree is $d=1$. According to \cite{Schroeer 2010},   Theorem 3.3 
the singular locus $\Sing(X)\subset X$ has codimension $d=1$. 
It follows that $\Sing(G)\subset G$ is not an equality.

Seeking a contradiction, we now assume that $Z=\Sing(G)$ has codimension $\geq 2$. Then the scheme $G$ is normal,
by Serre's Criterion. Choose a normal compactification $Y=\bar{G}$.
The canonical map $k\ra H^0(Y,\O_Y)$ is bijective, because we have the rational point $e\in G$.
According to \cite{Schroeer 2010}, Lemma 1.3 the base-change $Y\otimes_kk(\omega^{1/p})$ remains integral.
On the other hand, we just saw that $G\otimes k(\omega^{1/p})$ is non-reduced, contradiction.
\qed

\section{Witt algebras}
\mylabel{Witt algebras}

We keep the notation from the previous section. Our goal now is to understand
the restricted Lie algebra $\lieg=\lieg_\omega$, or equivalently the Frobenius kernel, attached to the 
automorphism group scheme $G=G_\omega$ of the spectrum of the ring $L=L_\omega=k[t]/(t^p-\omega)$. 

From $G=\Aut_{\Spec(L)/k}$ we get an  identification $\lieg=\Der_k(L)$. Any $k$-derivation  $\delta:L\ra L$ can be seen
as an $L$-linear map $\Omega^1_{L/k}\ra L_\omega$. The module of K\"ahler differentials
is a free  $L$-module of rank one, generated by $dt$. Let $\partial\in \lieg$ be the dual basis vector.
In turn we get the canonical $k$-basis $t^i\partial$, with $0\leq i\leq p-1$,
and the Lie bracket is given by 
$$
[t^i\partial,t^j\partial]=(j-i)t^{i+j-1}\partial.
$$
Using this relation with $i=0$, and also with $j=0$, together with $[t^{p-1}\partial,t\partial]=(2-p)t^{p-1}$, one easily sees
that $\lieg$ is simple, provided $p\neq 2$.

The $p$-map    $(f\partial)^\pmap=(f\partial)^p$ is the $p$-fold composition in $\End_k(L)$. It    can be made explicit as follows: 
For each truncated polynomial $f=\sum_{i=0}^{p-1}\lambda_it^i$, we write $f^{p-1}=\sum_{i=0}^{p-1} C_it^i$,
where the $C_i\in k$ 
are   certain   polynomial expressions in the  coefficients $\lambda_0,\ldots,\lambda_{p-1}$ and $\omega$, which  also depend on  the prime $p>0$.
Set  $C=C_{p-1}$. For example, with $p=5$ the polynomial $C$  becomes
\begin{gather*}
(\lambda_0^3\lambda_4 + 2\lambda_0^2\lambda_1\lambda_3 + \lambda_0^2\lambda_2^2 + 2\lambda_0\lambda_1^2\lambda_2 + \lambda_1^4) +  
\omega(2\lambda_0\lambda_1\lambda_4^2 +     4\lambda_0\lambda_2\lambda_3\lambda_4 +  4\lambda_0\lambda_3^3 + \\ 2\lambda_1^2\lambda_3\lambda_4 + 2\lambda_1\lambda_2^2\lambda_4 + 2\lambda_1\lambda_2\lambda_3^2 +     4\lambda_2^3\lambda_3)  +  
\omega^2 (4\lambda_2\lambda_4^3 + \lambda_3^2\lambda_4^2).
\end{gather*}

\begin{proposition}
\mylabel{description p-power}
We have 
$(f\partial)^\pmap = C \cdot f\partial$ for every element $f\partial\in\lieg$.  Moreover,
the factor $C$   is homogeneous of degree $p-1$ in the coefficients
of $f=\sum\lambda_it^i$.
\end{proposition}

\proof
Clearly we have $\partial^p=0$.
According to Hochschild's Formula (\cite{Hochschild 1955}, Lemma 1) the $p$-fold composition of $f\partial$ is given by 
$$
(f\partial)^p =  f^p\partial^p +  g\partial = g\partial,
$$
where $g=(f\partial)^{p-1}(f)$.  Consider the differential operator $D=\partial f\partial\ldots f\partial$,
where the number of $\partial$-factors is $p-1$, such that $g=fD(f)$.
According to \cite{Evans; Fuchs 2002}, Theorem 2 we have $D(f)=-\partial^{p-1}(f^{p-1})$. Note that this  result
is purely formal, and holds in any $\FF_p$-algebra   with a chosen element $f$ and some derivation $\partial$.
Clearly, $\partial^{p-1}(t^i)=0$ for $0\leq i\leq p-2$, whereas $\partial^{p-1}(t^{p-1})=(p-1)!=-1$.
Summing up, we have $D(f)= C$, hence  $g=fC=Cf$,  and the statement on the $p$-map follows.
From $(\sum\lambda_it^i)^{p-1}=\sum C_it^i$ one immediately sees that each $C_i=C_i(\lambda_0,\ldots,\lambda_{p-1})$
is homogeneous of degree $p-1$. 
\qed

\medskip
This has a remarkable consequence:

\begin{corollary}
\mylabel{witt p-closed}
Every vector in the restricted Lie algebra  $\lieg$ is $p$-closed.  
\end{corollary}

So each non-zero vector $f\partial\in\lieg$ defines a subgroup scheme $H\subset G$ of order $p$.
Note that the   additive vectors   might be viewed as    rational points on the hypersurface of degree $p-1$
defined by the homogeneous equation $C(\lambda_0,\ldots,\lambda_{p-1})=0$.
For the primes $p=2$ and $p=3$ we get $C(\lambda_0+\lambda_1)=\lambda_1$ and $C(\lambda_0,\lambda_1,\lambda_2)=\lambda_1^2-\lambda_0\lambda_2$,
respectively, which then reveals   the  structure of $\lieg$:

\begin{corollary}
\mylabel{small primes}
For $p=2$ we have $\lieg\simeq k\rtimes\liegl_1(k)$. For $p=3$ one gets $\lieg\simeq\liesl_2(k)$.
\end{corollary}

\proof
In the first case, one easily checks that the linear bijection 
$\lieg\ra k\rtimes\liegl_1(k)$ given by 
$(a+bt)\partial\mapsto (a,b)$ respects bracket  and $p$-map.
In the second case the linear bijection 
$$
\lieg\lra\liesl_2(k),\quad (a+bt+ct^2)\partial\longmapsto\begin{pmatrix}b&a\\-c&-b\end{pmatrix}
$$
likewise respects bracket and $p$-map.
\qed

\medskip
Consider the  \emph{adjoint representation}
$\Ad:G\ra\Aut_{\lieg/k}$, which sends each $g\in G(R)$ to the derivative of the conjugacy map $c_g$ given by $x\mapsto gxg^{-1}$.
From \cite{Waterhouse 1971}, Theorem in Section 5.2 we get:

\begin{proposition}
\mylabel{adjoint bijective}
For $p\geq 5$ the adjoint representation
$\Ad:G\ra\Aut_{\lieg/k}$ is an isomorphism of group schemes.
\end{proposition}

So for $p\geq 5$  our  $G$ can be seen as the automorphism group scheme for the  ring $L$,
the group scheme $G$, and the restricted Lie algebra $\lieg$.
Consequently, the three conditions in Proposition \ref{twisted form} are also equivalent to $\lieg_\omega\simeq\lieg_{\omega'}$.
For $p=2,3$ the adjoint representation $G\ra\Aut_{\lieg/k}$ is not bijective, according to Proposition \ref{twisted forms}.

\medskip
We now come to the crucial   result of this paper:

\begin{theorem}
\mylabel{subalgebras}
Suppose  that the scalar $\omega\in k$ is not a $p$-power, 
that $k^\times=k^{\times(p-1)}$, and that the Brauer group $\Br(k)$ contains no element of order two. Then each subalgebra  $\lieg'\subset\lieg$  
of dimension $1\leq n\leq p-1$ is  isomorphic to either  $k$ or  $\liegl_1(k)$ or $k\rtimes\liegl_1(k)$ or $\liesl_2(k)$.
\end{theorem}

\proof
In the special case $p=2$ the dimension of $\lieg'$ must be $n=1$, and it follows that $\lieg'$ is a twisted
form of $k$ or $\liegl_1(k)$. According to Proposition \ref{twisted forms}, all such twisted forms are trivial,
so our assertion indeed holds.

From now on we assume $p\geq 3$.
Recall that $\lieg=\Lie(G)$ is a twisted form of $\lieg_0=\Lie(G_0)$, where $G_0$ is the automorphism group
scheme for the spectrum of $L_0=k[t]/(t^p)$.  Let $\lieg_{0,\red}$ be the subalgebra corresponding to 
the reduced part $G_{0,\red}$.  According to Proposition \ref{twisting adjoint} below, there is no vector $x\neq 0$ in $\lieg$
such that $x\otimes 1\in \lieg\otimes k(\omega^{1/p})$ is contained in $\lieg_{0,\red}\otimes k(\omega^{1/p})$.
In particular, the latter does not contain  the base-change $\lieg'\otimes k(\omega^{1/p})$.

It follows that the further base-change $\lieg'\otimes k^\alg$ is not contained in $\lieg_{0,\red}\otimes k^\alg$.
Such subalgebras were studied by Premet and Stewart in \cite{Premet; Stewart 2019},  Section 2.2:
They remark on page 971 that a subalgebra  in $\lieg_0\otimes k^\alg$ is not contained in $\lieg_{0,\red}\otimes k^\alg$
if and only if it does not preserve any proper non-zero ideal, and call such subalgebras \emph{transitive}.
We thus may apply loc.\ cit., Lemma 2.2 and infer that
$\lieg'$ is a twisted form of $k$ or  $\liegl_1(k)$ or $k\rtimes\liegl_1(k)$ or $\liesl_2(k)$.
By assumption, the groups $k^\times/k^{\times(p-1)}$ and $\Br(k)[2]$ vanish.
According to Lemma \ref{twisted forms}, the four restricted Lie algebras in question have no twisted forms over our field $k$,
thus $\lieg'$ is isomorphic to one of them.
\qed


\medskip
Consequently, for every $\lieg'\subset\lieg$ as above over any ground field $k$, one finds a finite separable extension 
so that the base-change of $\lieg'$ belongs to the given list.

\section{Twisting adjoint representations}
\mylabel{Twisting}

%
%
%
We keep the assumption of the preceding section, and establish the crucial ingredients  for the proof of Theorem \ref{subalgebras}.
Recall  that we are in characteristic $p>0$, and
that  $G=G_\omega$ is the automorphism group scheme of the spectrum of $L=L_\omega=k[t]/(t^p-\omega)$,
for some scalar $\omega\in k$.
The resulting  restricted Lie algebra $\lieg=\Lie(G)$ is the $p$-dimensional vector space $\Der_k(L)$, which comprises
the derivations $f(t)\partial$, where $f=\sum_{i=0}^{p-1}\mu_it^i$ is a truncated polynomial.
Moreover,   the group elements $g\in G(R)$ act from the left on  the spectrum of $L\otimes_k R$, and from the right
on the coordinate ring  $L\otimes_kR$  via the substitution $t\mapsto \varphi_g(t)$,
for the corresponding truncated polynomial $\varphi_g(t)=\sum_{i=0}^{p-1}\lambda_it^i$. The  coefficients define an
embedding $G\subset\AA^p$ of the underlying scheme.  For each $g\in G(R)$,   write $c_g$ for the induced inner automorphism
$x\mapsto gxg^{-1}$. The resulting  conjugacy map $c:G\ra\Aut_{G/k}$ is given in terms of truncated polynomials by the formula
$$
\varphi_{gxg^{-1}}(t) = \varphi_{g^{-1}}(\varphi_x(\varphi_g(t))).
$$
By functoriality, the elements $c_g\in\Aut_{G/k}(R)$ induce an automorphism $\Ad_g=\Lie(c_g)$ of  $\lieg\otimes_kR$,
which defines the adjoint representation $\Ad:G\ra \Aut_{\lieg/k}$.

\begin{proposition}
\mylabel{adjoint representation}
Let  $g\in G(R)$, and write  $\varphi(t)=\varphi_{g^{-1}}(t)$ for the truncated polynomial of the inverse $g^{-1}$. Then the formal derivative
$\varphi'(t)$ is a unit in the ring $L\otimes_kR$, and  for each $f(t)\partial\in\lieg\otimes_kR$ we have
$$
\Ad_g(f(t) \partial)=\ \frac{f(\varphi(t))}{\varphi'(t) }\partial.
$$
\end{proposition}

\proof
By definition, the element $f(t)\partial\in\lieg\otimes_kR\subset G(R[\epsilon])$ acts on the algebra $L\otimes R[\epsilon]$ via
$$h(t)\longmapsto h(t)+\epsilon f(t)\partial(h) = h(t)+\epsilon f(t)h'(t).
$$
Thus the adjoint $\Ad_g(f(t)\partial)=g^{-1}\circ f(t)\partial\circ g$ is given by the following composition:
\begin{equation}
\label{derived conjugacy}
t\longmapsto \varphi_g(t)\longmapsto \varphi_g(t) + \epsilon f(t)\varphi_g'(t)\longmapsto
\varphi_g(\varphi_{g^{-1}}(t)) + \epsilon f(\varphi_{g^{-1}}(t))\varphi_g'(\varphi_{g^{-1}}(t)).
\end{equation}
For a moment, let us regard the truncated polynomials $\varphi_g(t)$ and $\varphi_{g^{-1}}(t)$ as elements in the polynomial ring $R[t]$.
Then $t=\varphi_g(\varphi_{g^{-1}}(t)) + (t^p-\omega)h(t)$ for some polynomial $h(t)$. Taking formal derivatives and applying the chain rule,
we obtain 
$$
1=\varphi_g'(\varphi_{g^{-1}}(t))\cdot \varphi'_{g^{-1}}(t) + (t^p-\omega)h'(t).
$$
This gives $1=\varphi_g'(\varphi_{g^{-1}}(t))\cdot \varphi'_{g^{-1}}(t)$ in the truncated polynomial ring  $L\otimes_kR=R[t]/(t^p)$.
It follows that $\varphi(t)=\varphi'_{g^{-1}}(t)$ is a unit, with inverse $\varphi_g'(\varphi_{g^{-1}}(t))$.
Substituting for the term on the right in \eqref{derived conjugacy} gives the desired formula for $\Ad_g(f(t) \partial)$.
\qed

\medskip
Now consider the additive vector $\partial\in\lieg$, which corresponds
to an inclusion of the infinitesimal group scheme $H=\alpha_p$ into the group scheme $G$.
The $R$-valued points $h\in H(R)=\{\lambda\in R\mid\lambda^p=0\}$ correspond to truncated polynomials $\varphi_h(t)=t+\lambda$.
The inverse $R$-valued point  has $\varphi_{h^{-1}}=t-\lambda$,
with formal derivative $\varphi_{h^{-1}}'(t)=1$. This immediately gives:

\begin{corollary}
\mylabel{additive adjoint}
With the above notation, we have $\Ad_h(f(t)\partial)= f(t-\lambda)\partial$ for every element $f(t)\partial\in \lieg\otimes_kR$.
\end{corollary}

Recall that $G=G_\omega$ depends on some scalar $\omega\in k$, and is a twisted form of $G_0$. The latter coincides with
its own Frobenius pullback.
By Proposition \ref{reduced part}, the reduced part $G_{0,\red}$ is a non-normal subgroup scheme.
Recall that the embedding $G_0\subset \AA^p$ is given by $\lambda_0^p=0$ and $\lambda_1\neq 0$,
such that $G_{0,\red}$ is defined by $\lambda_0=0$ and $\lambda_1\neq 0$.
Write $\lieg_{0,\red}\subset\lieg_0$ for the resulting subalgebra, which comprises the derivations $f\partial$
where the truncated polynomial $f=\sum_{i=0}^{p-1}\lambda_it^i$ has $\lambda_0=0$.
Write $H_0\subset G_0$ for the copy of $\alpha_p$ given by the additive vector $\partial\in\lieg_0$.
 
Now suppose that our ground field $k$ is imperfect,   that our scalar $\omega\in k$   is not a $p$-power, 
and consider the resulting field extension $k(\omega^{1/p})$.   In light of Lemma \ref{torsor description}, we may endow its spectrum $T$
with the structure of an $H_0$-torsor.  Lemma \ref{twisted automorphisms} gives an identification ${}^T\lieg_0=\lieg_\omega$, and thus 
an identification $\lieg_0\otimes_kk(\omega^{1/p})=\lieg_\omega\otimes_kk(\omega^{1/p})$
The following fact was a crucial ingredient for the proof of Theorem \ref{subalgebras}:

\begin{proposition}
\mylabel{twisting adjoint}
The twisted form $\lieg_\omega$ contains no   vector $x\neq 0$ such that the induced vector $x\otimes 1$
inside $\lieg_0\otimes_kk(\omega^{1/p})=\lieg_\omega\otimes_kk(\omega^{1/p})$ is contained in the base-change  $\lieg_{0,\red}\otimes_kk(\omega^{1/p})$.
\end{proposition}

\proof
Setting $V=\lieg_0$ and $V'=\lieg_{0,\red}$, we see that that action of $H_0=\alpha_p$ via the adjoint representation  $G_0\ra \Aut_{\lieg_0/k}$
is exactly as described in Proposition \ref{no subspace}, and the assertion follows.
\qed

\section{Subalgebras}
\mylabel{Subalgebras}

Throughout this section, $k$ is  a field of characteristic $p>0$, and $k\subset E$ is a field extension.
Suppose we have a group scheme $H$ of finite type over $k$,   a group scheme $G$
of finite type over $E$, and a  homomorphism $f:H\otimes_kE\ra G$.
We shall see that in certain circumstances,  important structural properties of the Frobenius kernel $G[F]$ are inherited to $H[F]$.

Consider the finite-dimensional restricted Lie algebra $\lieh=\Lie(H)$ over $k$
and $\lieg=\Lie(G)$ over $E$.
Our homomorphism   of group schemes  induces an $E$-linear homomorphism  
$$
\Lie(f):\lieh\otimes_kE\ra\lieg,\quad x\otimes\alpha\longmapsto \alpha x,
$$
of restricted Lie algebras
which  corresponds to  a $k$-linear homomorphism 
 $\lieh\ra\lieg$ of restricted Lie algebras.
We are mainly interested in the case that $E$ is the function field of an integral $k$-scheme $X$ of finite type,
such that   $\lieg$ is an infinite-dimensional $k$-vector space.
Set $N=\Kernel(f)$, with Lie algebra $\lien=\Kernel(\Lie(f))$.

\begin{proposition}
\mylabel{characterization injective}
The following are equivalent:
\begin{enumerate}
\item The $k$-linear homomorphism $\lieh\ra\lieg$ is injective.
\item For every  non-trivial subgroup scheme $H'\subset H$ that is minimal with respect to inclusion,
the base-change $H'\otimes_kE$ is not contained in the kernel  $N\subset H\otimes_kE$.
\end{enumerate}
\end{proposition}

\proof
We prove the contrapositive: Suppose $\lieh\ra\lieg$ is not injective.
Inside the kernel, choose a subalgebra $\lieh'\neq 0$ that is minimal with respect to inclusion.
Then the induced $E$-linear map $\lieh'\otimes_kE\ra\lieg$ is zero.  
By the Demazure--Gabriel Correspondence, the corresponding subgroup scheme $H'\subset H$
of height one is minimal with respect to inclusion, and $H'\otimes_kE\subset N$.
Conversely, suppose $H'\otimes_kE\subset N$ for some $H'$ as in (ii). Choose some non-zero vector $x$ from $\lieh'=\Lie(H')$.
By construction, it lies in the kernel for $\lieh\ra\lieg$.
\qed

\medskip
We now suppose that the above equivalent conditions hold, and regard the injective map as an inclusion $\lieh\subset\lieg$.
To simplify exposition, we also assume that $k$ is algebraically closed, and that $\lieh$   contains an $E$-basis for $\lieg$. In other words,
the induced linear map $\lieh\otimes_kE\ra \lieg$ is surjective. Note  
that this $E$-linear map is usually   \emph{not injective}.
However, we shall see that  important structural properties of $\lieg$ transfer to $\lieh$.
We start with a series of three elementary but useful observations:

\begin{lemma}
\mylabel{all p-closed}
If every vector in $\lieg$ is $p$-closed, the same holds for every vector in $\lieh$.
\end{lemma}

\proof
Fix some non-zero $x\in\lieh$. By assumption we have $x^\pmap=\alpha x$ for some   $\alpha\in E$, 
and our task is to verify that this scalar already lies in $k$.
Since the latter is algebraically closed, it is enough to verify that $\alpha$ is algebraic over $k$.
By induction on $i\geq 0$ we get $x^{[p^i]}=\alpha^{n_i}x$ for some   strictly increasing sequence   $0=n_0<n_1<\ldots$
of integers.
Since $\dim_k(\lieh)<\infty $ there is a non-trivial relation $\sum_{i=0}^r\lambda_i x^{[p^i]}$ for some $r\geq 0$
and some coefficients $\lambda_i\in k$. This gives $\sum  \lambda_i\alpha^{n_i}x=0$.
Since $x\neq 0$ we must have $\sum\lambda_i\alpha^{n_i}=0$,
hence $\alpha\in E$ is algebraic over $k$.
\qed

\begin{lemma}
\mylabel{toral rank}
The restricted Lie algebras $\lieg$ and $\lieh$ have the same toral rank, and the kernel $\lien$ for 
  $\lieh\otimes_kE\ra \lieg$ has toral rank $\rho_t(\lien)=0$.
\end{lemma}

\proof
It follows from  \cite{Block; Wilson 1988}, Lemma 1.7.2 that  $\rho_t(\lieh)=\rho_t(\lieg)+\rho_t(\lien)$,
and in particular  $\rho_t(\lieg)\leq \rho_t(\lieh)$.
For the reverse inequality, suppose there are  $k$-linearly independent vectors $x_1,\ldots,x_r\in \lieh$ with $[x_i,x_j]=0$ and $x_i^\pmap=x_i$.
We have to check that the vectors are $E$-linearly independent. Suppose there is a non-trivial relation.
Without loss of generality, we may assume that $x_1,\ldots,x_{r-1}$ are $E$-linearly independent,
and that $x_r=\sum_{i=1}^{r-1}\lambda_ix_i$ for some coefficients $\lambda_i\in E$.
From the axioms of the  $p$-map  we get
$$
\sum \lambda_ix_i=x_r=x_r^\pmap=(\sum\lambda_ix_i)^\pmap= \sum\lambda_i^px_i^\pmap = \sum\lambda_i^px_i.
$$
Comparing coefficients gives $\lambda_i^p = \lambda_i$. Thus $\lambda_i$ lie in the prime field, in particular in $k$.
In turn, the vectors are $k$-linearly dependent, contradiction. 
\qed

\medskip
Let us call the restricted Lie algebra  $\lieh$   \emph{simple} if it is non-zero, and contains no ideal besides
$\liea=0$ and $\liea=\lieh$. 
 
\begin{lemma}
\mylabel{simple algebra}
Suppose there is a restricted Lie algebra $\lieh'$ over $k$ such that
$\lieg$ is a twisted form of the base-change $\lieh'\otimes_kE$.
If $\lieh'$ is simple of dimension $n'\geq 2$, we must have $\lieh\simeq\lieh'$ and $\lieg\simeq\lieh'\otimes_kE$.
\end{lemma}

\proof
Let $H$ and $H'$ be the finite group schemes of height one corresponding to the restricted Lie algebras $\lieh$ and $\lieh'$,
respectively. Consider the Hom scheme $X\subset \Hilb_{H\times H'}$ of surjective homomorphisms $H\ra H'$.
By assumption, this scheme contains a point with values in the algebraic closure $E^\alg$.
By Hilbert's Nullstellensatz, there must be a point with values in $k$, hence 
there there is a surjective homomorphism $H\ra H'$. It corresponds to a  short exact sequence
of restricted Lie algebras
$$
0\lra \liea\lra \lieh\lra\lieh'\lra 0.
$$
We claim that the  ideal $\ideala$ vanishes. Suppose this is not the case.
Clearly, $\lieh'$ and $\lieg$ have the same toral rank. By Lemma \ref{toral rank}, also $\lieg$ and $\lieh$ have
the same toral rank. According to \cite{Block; Wilson 1988}, Lemma 1.7.2 we have $\rho_t(\liea)=0$,
so the $p$-map on $\liea$ is nilpotent.
On the other hand, the $p$-map on $\lieh'$ is not nilpotent, because the Lie algebra is simple
of dimension $\dim(\lieh')\geq 2$. The    same holds for $\lieg$,
and we infer that the  induced map $\liea\otimes_kE\ra\lieg$ is not surjective.
Its image $\lieb\subsetneq\lieg$ is non-zero, because   $\lieh\subset\lieg$. 
Since $\lieg$ is simple, there   are elements $x\in \lieb$ and $y\in \lieg$
with $[x,y]\not\in \lieb$. Such vectors may be chosen with  $x\in \liea$ and  $y\in \lieh$, because  
$\liea\subset \lieb$ and $\lieh\subset \lieg$ contain $E$-bases.
Consequently, $\liea\subset\lieh$ is not an ideal, contradiction.

This shows that $\lieh=\lieh'$. 
In particular, $\lieh$ and $\lieg$ have the same vector space dimension, so our surjection 
$\lieh\otimes_kE\ra \lieg$ must be bijective. Our assertions follow.
\qed

\medskip
For each  $a\in\lieh$, the Lie bracket $\ad_a(x)=[a,x]$ defines a $k$-linear endomorphism of $\lieh$,
but also an $E$-linear endomorphism of $\lieg$. Write $\ad_{\lieh,a}$ and $\ad_{\lieg,a}$ for the respective maps,
and  $\mu_{\lieh,a}(t)\in k[t]$ and $\mu_{\lieg,a}(t)\in E[t]$ for the resulting  minimal polynomials.

\begin{lemma}
\mylabel{minimal polynomial}
We have $\mu_{\lieh,a}(t)=\mu_{\lieg,a}(t)$. In particular, the endomorphism
$\ad_{\lieg,a}$   is trigonalizable, and its eigenvalues  coincide with those
of $\ad_{\lieh,a}$. Moreover, the former is diagonalizable if and only if this holds for the latter.
\end{lemma}

\proof
The surjection $\lieh\otimes_kE\ra\lieg$ already reveals that $\mu_{\lieg,a}(t)$ divides $\mu_{\lieh,a}(t)$.
The latter decomposes into linear factors over $k$, because this field is algebraically closed.
We conclude that $\mu_{\lieg,a}(t)=\sum\lambda_it^i$ actually lies in $k[t]$, and decomposes into linear factors over $k$.
Moreover, for each vector $x$ from $\lieh\subset\lieg$ we have  $\sum\lambda_i\ad_{\lieh,a}^i(x)=0$, hence $\mu_{\lieh,a}(t)$ divides
$\mu_{\lieg,a}(t)$. In turn, the two minimal polynomials coincide.
The remaining assertions follow immediately.
\qed

\medskip
We now consider some   special cases for  $\lieg$, and deduce  structure results for $\lieh$.
Recall that $k^n$ denotes the $n$-dimensional restricted Lie algebra over $k$ with trivial bracket 
and $p$-map. The following fact is  obvious:

\begin{proposition}
\mylabel{additive}
If $\lieg$ is isomorphic to $E^m$ then the restricted Lie algebra   $\lieh$ is isomorphic to   $k^n$
for some integer $n\geq m$.
\end{proposition}
 
Recall that $k^n\rtimes_\varphi\liegl_1(k)$ denotes the semidirect product formed with respect to 
the homomorphism $\varphi:\liegl_1(k)\ra\liegl(k^n)=\Der'_k(k^n)$ that sends scalars to scalar matrices.

\begin{proposition}
\mylabel{semidirect product}
If  $\lieg$ is isomorphic to  $E^m\rtimes\liegl_1(E)$, then  the restricted Lie algebra $\lieh$ is isomorphic to $ k^n\rtimes\liegl_1(k)$
for some $n\geq m$.
\end{proposition}

\proof
Without loss of generality we may assume $\lieg=E^m\rtimes\liegl_1(E)$. First recall that    bracket and $p$-map are given by the formulas
\begin{equation}
\label{semidirect formulas}
[v+\lambda e,v'+\lambda'e]=\lambda v'-\lambda'v\quadand (v+\lambda e)^\pmap = \lambda^{p-1}(v+\lambda e),
\end{equation}
where $v\in E^n$, and  $e\in \liegl_1(E)$ denotes the unit.
In particular, each vector is $p$-closed. Moreover, $a=v+\lambda e$ is multiplicative if and only if $\lambda\neq 0$,
and $a^\pmap=a$ if and only if  $\lambda\in\mu_{p-1}(E)$. 
For any such vector,  we see that the endomorphism $\ad_a(x)=[a,x]$ is diagonalizable,
and $E^n\subset\lieg$ is the eigenspace with respect to the eigenvalue $\alpha=\lambda$, whereas the line $Ea\subset\lieg$ is
the eigenspace for $\alpha=0$. 

From the extension $0\ra E^m\ra\lieg\ra\liegl_1(E)\ra 0$ one   sees that $\lieg$ has toral rank one.
According to Lemma \ref{toral rank}, the same holds for $\lieh$.
Choose some non-zero vector $a\in \lieh$ with $a^\pmap=a$. Then $a=v+\lambda e$ for some $\lambda\in\mu_{p-1}(E)\subset k^\times$.
Replacing $a$ by $\lambda^{-1}a$ we may assume $\lambda=1$.
By Lemma \ref{minimal polynomial}, the adjoint representation $\ad_{\lieh,a}$   is diagonalizable, with eigenvalues $\alpha=0$ and $\alpha=1$.
Let $\lieh=U_0\oplus U_1$ be the corresponding eigenspace decomposition.
Then $U_0$ lies  in the corresponding eigenspace for $\ad_{\lieg,a}$, which is $E^n\subset\lieg$. 
It follows that $U_0$ has trivial Lie bracket and $p$-map. The choice of a $k$-basis gives $U_0=k^n$ for some $n\geq 0$.
Likewise, $U_1$ is contained in $Ea$. Thus the bracket vanishes on $U_1$, and the $p$-map is injective.
Using Lemma \ref{toral rank}, we infer that $U_1=ka$. The vector space decomposition $\lieh=U_0\oplus U_1$ thus
becomes a semidirect product $\lieh=k^n\rtimes\liegl_1(k)$. We must have $m\geq n$ because the map $\lieh\otimes_kE\ra\lieg$ is surjective.
\qed

\medskip
Recall that  $\liesl_2(E)$ is simple for $p\geq 3$. Using Lemma \ref{simple algebra}, we immediately obtain:

\begin{proposition}
\mylabel{sl}
Suppose $p\geq 3$. If $\lieg$ is isomorphic to a twisted form of $\liesl_2(E)$, then the restricted Lie algebra   $\lieh$ is isomorphic to $\liesl_2(k)$. 
\end{proposition}

\section{Structure results for Frobenius kernels}
\mylabel{Structure results}

We now come to our main result.
Let $k$ be an algebraically closed field of characteristic $p>0$,
and $X$ be a proper integral scheme or more generally a proper integral algebraic space, 
$H=\Aut_{X/k}[F]$ be the Frobenius kernel for the automorphism group scheme,
and $\lieh=H^0(X,\Theta_{X/k})$ the corresponding restricted Lie algebra over $k$.
Let $H_F=H\otimes_kF$ be the base-change to the function field $F=k(X)$,  
and $H_F^\inert\subset H_F$ the inertia subgroup scheme for the rational point
in the spectrum of $F\otimes_E F$, with corresponding restricted Lie algebra $\lieh_F^\inert\subset\lieh_F$.
Recall that the foliation rank $r\geq 0$ is given by
$$
r=\dim(\lieh_F/\lieh_F^\inert)=\dim(\Omega^1_{F/E})\quadand [H_F:H_F^\inert]= [F:E] = p^r,
$$
where $E=F^\lieh$ is the kernel for all derivations $D:F\ra F$ from
the Lie algebra $\lieh$. This is nothing but the function field 
$E=k(Y)$  of the quotient $Y=X/H$.

\begin{theorem}
\mylabel{structure result}
Suppose that the proper integral scheme $X$ has foliation rank $r\leq 1$.  Then the Frobenius kernel $H=\Aut_{X/k}[F]$ is isomorphic to the Frobenius
kernel  of one of the following three basic types of group schemes:
$$
\SL_2\quadand \GG_a^{\oplus n} \quadand \GG_a^{\oplus n}\rtimes\GG_m,
$$
for some integer $n\geq 0$.
\end{theorem}

\proof
The case $r=0$ is trivial, so we assume $r=1$, such that $\lieh\neq 0$.
By assumption the subfield $E=F^\lieh$ has $[F:E]=p$, and we thus have $F=E[T]/(T^p-\omega)$ for some
$\omega\in E$. Thus the restricted Lie algebra $\lieg=\Der_E(F)$ is a twisted form of the Witt algebra $\lieg_0$ over $E$.
By construction, we have an inclusion $\lieh\subset\lieg$. 

Suppose first that the induced homomorphism $\lieh\otimes E\ra \lieg$ is surjective, such that the results of Section \ref{Subalgebras}
apply. Suppose first $p\leq 3$. Then $\lieg$ is isomorphic to either $E\rtimes\liegl_1(E)$ or $\liesl_2(E)$, by Corollary  \ref{small primes}, 
and our assertion follows from Proposition \ref{semidirect product} and Proposition \ref{sl}.
The case $p\geq 5$ actually does not occur: Then the Witt algebra $\lieh'=\Der_k(k[t]/(t^p))$ is simple, 
as   remarked at the beginning of Section \ref{Witt algebras},
and $\lieg$ is a twisted form of $\lieh'\otimes _kE$. It follows from Lemma  \ref{simple algebra} that $\lieg\simeq \lieh'\otimes_kE$.
Combining Proposition \ref{adjoint bijective} and Proposition \ref{twisted form}, we get $E[T]/(T^p-\omega)\simeq E[T]/(T^p)$, contradiction.

It remains to treat the case that $\lieh\otimes E\ra \lieg$ is not surjective. Then $\lieg'=\lieh\cdot E$
is a restricted subalgebra of dimension $1\leq n\leq p-1$. 
We now replace the field $E$ by some separable closure $E^\sep$, and likewise $F$ by $F\otimes_EE^\sep$.
According to Theorem \ref{subalgebras} the restricted Lie algebra $\lieg'$ is isomorphic to either
$\liesl_2(E)$ or $E$ or $\liegl_1(E)$ or $E\rtimes\liegl_1(E)$.
By the results in Section \ref{Subalgebras}, this ensures that
$\lieh$ is isomorphic to $\liesl_2(k)$ or $k^n$ or $k^n\rtimes\liegl_1(k)$    for some $n\geq 0$.
These are the Frobenius kernels for the group schemes in question, and our assertion follows.
\qed

\medskip
In the former case, the Frobenius kernel is $\liesl_2(k)$. This indeed occurs for $X=\PP^1$.
In the latter two cases, the respective Frobenius kernels are $\alpha_p^{\oplus n}$ and   $\alpha_p^{\oplus n}\rtimes\mu_p$.
With Proposition \ref{normal surface}, we immediately get the following consequence:
 
\begin{corollary}
\mylabel{structure for surfaces}
Suppose that $X$ is a proper normal surface with  $h^0(\omega^{\vee}_X)=0$.
Then  $H=\Aut_{X/k}[F]$ is isomorphic to the Frobenius
kernel  of one of the   three basic types of group schemes in the Theorem.
\end{corollary}

This applies in particular to smooth surface $S$  of Kodaira dimension $\kod(S)\geq 1$,
to surfaces of general type   and their minimal models,
or normal  surfaces $X$ with $c_1=0$ and $\omega_X\neq\O_X$ having at most rational double points.

\section{Canonically polarized surfaces}
\mylabel{Canonically polarized}

Let $k$ be a ground field of characteristic $p>0$,
and $X$ be a proper normal surface with $h^0(\O_X)=1$ whose complete local rings $\O_{X,a}^\wedge$ are complete intersections.
Then the cotangent complex $L^\bullet_{X/k}$ is perfect, and we obtain two \emph{Chern numbers}
$$
c_1^2=c_1^2(L^\bullet_{X/k}) \quadand c_2=c_2(L^\bullet_{X/k}),
$$
as explained by Ekedahl, Hyland and Shepherd--Barron \cite{Ekedahl; Hyland; Shepherd-Barron 2012}, Section 3.
In some sense, these integers are the most fundamental numerical invariants of the surface $X$.
Note that $c_1^2$ is nothing but the self-intersection number $K_X^2=(\omega_X\cdot\omega_X)$ of the dualizing sheaf.
If the singularities are also rational, hence rational double points, the Chern numbers of $X$ coincide
with the Chern numbers of the minimal resolution of singularities $S$, according to loc.\ cit.\ Proposition 3.12 and Corollary 3.13.
For more details on rational double points, we refer to
\cite{Lipman 1969} and \cite{Artin 1977}.

Recall that a \emph{canonically polarized surface} is the canonical model $X$
of a smooth surface $S$ of general type. Then $\omega_X$ is ample,  all local rings $\O_{X,a}$ are either regular 
or rational double points,  and the above applies. Let us record the following facts:

\begin{lemma}
\mylabel{inequalities}
Suppose that $X$ is canonically polarized. Then  
$$
\chi(\O_X)=\frac{1}{12}(c_1^2+c_2)\quadand h^0(\omega_X)\leq \frac{1}{2}(c_1^2+4)\quadand c_2\leq 5c_1^2 + 36.
$$
\end{lemma}

\proof
Let $f:S\ra X$ be the minimal resolution of singularities. Then $S$ is a smooth minimal surface of general type, 
and the first formula holds for $S$ instead of $X$ by Hirzebruch--Riemann--Roch.
We already observed that the surfaces $X$ and $S$ have the same Chern numbers,
and the structure sheaves have the same cohomology. Thus the formula also holds for $X$.
In particular we have $h^0(\omega_X)=h^2(\O_X) = h^2(\O_S)=h^0(\omega_S)$, and Noether's Inequality (for example \cite{Liedtke 2013b}, Section 8.3) for 
the minimal surface of general type $S$ gives
the second formula. This ensures
$$
\chi(\O_X) = 1-h^1(\omega_X) + h^0(\omega_X)\leq 1 + h^0(\omega_X) = 1 + h^0(\omega_S)\leq  (c_1^2+6)/2.
$$
Combining with $\chi(\O_X)=(c_1^2+c_2)/12$  we get the third inequality.
\qed
%
%
%

\begin{theorem}
\mylabel{bound result}
Let $X$ be canonically polarized surface, with Chern numbers $c_1^2, c_2$.
Then  the Lie algebra $\lieh=H^0(X,\Theta_{X/k})$ for the  Frobenius kernel $H=\Aut_{X/k}[F]$ has the property
$\dim(\lieh)\leq \Phi(c_1^2,c_2)$ for the polynomial
$$
\Phi(x,y) =
\begin{cases}
\frac{1}{144}(73x+y)^2-1	& \text{if $c_1^2\geq 2$;}\\
\frac{1}{144}(121x+y)^2-1	& \text{if $c_1^2=1$.}
\end{cases}
$$
Moreover, we also have the weaker bound $\dim(\lieh)\leq \Psi(c_1^2)$ with the polynomial 
$$
\Psi(x)=
\begin{cases}
\frac{169}{4}x^2+39x+8		& \text{if $c_1^2\geq 2$;}\\
\frac{441}{4}x^2 + 63x+8	& \text{if $c_1^2=1$.}
\end{cases}
$$
\end{theorem}

\proof
Fix some $m\geq 3$.  Serre Duality gives $h^i(\omega_X^{\otimes m})=h^{2-i}(\omega_X^{\otimes 1-m})$.
This vanishes for $i=2$, because $\omega_X$ is ample, and also for $i=1$ by \cite{Ekedahl 1988}, Chapter II, Theorem 1.7.
Thus we have $h^0(\omega_X^{\otimes m})=\chi(\omega_X^{\otimes m})$, and Riemann--Roch gives
\begin{equation}
\label{ekedahl bound}
h^0(\omega_X^{\otimes m})=  \chi(\O_X) + \frac{1}{2}(m^2-m) c_1^2.
\end{equation} 
According to \cite{Ekedahl 1988}, Chapter III, Theorem 1.20 the invertible sheaf $\omega_X^{\otimes m} $ is very ample for  $c_1^2\geq 2$ and $m=4$,
or   $c_1^2=1$ and $m=5$. It then defines a closed embedding
$X\subset \PP^n$ with $\omega_X=\O_X(1)$ and  $n+1=h^0(\omega^{\otimes m}_X)$.

The following argument, which give a better estimate than our original reasoning, was kindly communicated by the referee:
The canonical linearization of $\omega_X$ and its power $\omega_X^{\otimes m}$ yields a homomorphism $\Aut_{X/k}\ra\Aut_{\PP^n/k}=\PGL_{n+1,k}$
such that the inclusion $X\subset \PP^n$ is equivariant with respect to the action  of  $G=\Aut_{X/k}$.
In particular, the homomorphism of group schemes is a closed embedding, and the  resulting inclusion of tangent spaces
$H^0(X,\Theta_{X/k})\subset H^0(\PP^n,\Theta_{\PP^n/k})$ gives the estimate $h^0(\Theta_{X/k})\leq h^0(\omega_X^{\otimes m})^2-1$.
Substituting \eqref{ekedahl bound} and using 
$\chi(\O_X)  =(c_1^2+c_2)/12$, we get
$$
h^0(\Theta_X) \leq \left( \frac{c_1^2+c_2}{12} +\frac{(m^2-m)}{2}c_1^2\right)^2-1=\frac{1}{144}\left((6m^2-6m+1)c_1^2+c_2\right)^2-1.
$$
Setting $m=4$ and $m=5$ we get the desired bound 
$\dim(\lieh)\leq \Phi(c_1^2,c_2)$. Finally,  the inequality  $c_2\leq 5c_1^2+36$ from  Lemma \ref{inequalities}  yields the weaker bound $\dim(\lieh)\leq \Psi(c_1^2)$.
\qed

\section{Examples}
\mylabel{Examples}

Let $k$ be a ground field of characteristic $p>0$.
In this section we give examples of canonically polarized surfaces $X$ where
the Frobenius kernel of the automorphism group scheme is isomorphic to $\alpha_p^{\oplus n}\rtimes\mu_p$ and $\alpha_p^{\oplus m}$.
Note that we do not have  examples where $\SL_2[F]$ occurs.

To start with,  view $\PP^2$ as the homogeneous spectrum of $k[T_0,T_1,T_2]$. Fix some $d\geq 1$,  set  $\shL=\O_{\PP^2}(d)$ and  consider the section
\begin{equation}
\label{defining relation}
s=T_0T_1T_2^{pd-2} +   T_1T_2T_0^{pd-2} + T_2T_0T_1^{pd-2} \in \Gamma(\PP^2,\shL^{\otimes p}).
\end{equation}
Regarded as     $ \shL^{\otimes-p}\ra\O_{\PP^2}$, this
endows the coherent sheaf $\shA=\bigoplus_{i=0}^{p-1}\mathscr{L}^{\otimes-i} $ with the structure of a  $\ZZ/p\ZZ$-graded $\O_{\PP^2}$-algebra,
and we define $X=\Spec(\shA)$ as the relative spectrum.
 
\begin{proposition}
\mylabel{example 1}
In the above setting, suppose  $p\neq 3$ and $d\geq 4$.  Then
$$
\lieh=k^n\rtimes \lieg\mathfrak{l}_1(k) \quadand \mathrm{Aut}_{X/k}[F] =\alpha_p^{\oplus n}\rtimes \mu_p,
$$
where $n=(d+1)(d+2)/2$. Moreover, $X$ is a canonically polarized surface 
with  Chern invariants   $c_1^2=p(pd-d-3)^2$ and $c_2=3p+dp(p-1)(pd-3)$.
\end{proposition}

\proof
Being locally a hypersurface in affine three-space, the scheme  $X$ is Gorenstein. According to \cite{Ekedahl 1988}, Chapter I, Proposition 1.7
the dualizing sheaf is given by $\omega_X=\pi^*(\omega_{\PP^2}\otimes \mathscr{L}^{p-1})$, which equals the pullback
of $\O_{\PP^2}(pd-d-3)$. The statement on $c_1^2$ follows. Using $d(p-1)-3\geq d-3\geq 1$ we see that $\omega_X$ is ample.
Since $\pi:X\ra\PP^2$ is finite, the Euler characteristic $\chi(\O_X)$ equals
$$
\sum_{i=0}^{p-1}\chi(\O_{\PP^2}(-id))=\sum_{i=0}^{p-1}\binom{2-id}{2}  = \frac{12p - 9d(p-1)p +d^2(p-1)p(2p-1)}{12}.
$$
Now suppose for the moment that we already know that $X$ is geometrically normal, with  only rational double points.
Then $X$ is a canonically polarized surface, and   Lemma \ref{inequalities} yields  the statement on $c_2$.

We proceed by  computing $\lieh=H^0(X,\Theta_{X/k})$ as a vector space.
The grading of the structure sheaf $\shA=\bigoplus_{i=0}^{p-1}\mathscr{L}^{\otimes-i} $ corresponds to an action of $G=\mu_p$
on the scheme $X$, with quotient $\PP^2$. Let $D:\O_X\ra \Theta_{X/k}$ be the corresponding multiplicative vector field,
and $\O_X(\Delta)\subset \Theta_{X/k}$ the saturation of the image, for some effective Weil divisor $\Delta\subset X$.
Lemma \ref{four term sequence} below gives  an exact sequence
$$
0\lra\O_X(\Delta)\stackrel{D}{\lra} \Theta_{X/k} \lra  \omega_X^{\otimes -1}(-\Delta).
$$
The term on the right has no non-zero global sections, because $\omega_X$ is ample, 
and consequently $H^0(X,\O_X(\Delta))=H^0(X,\Theta_{X/k})$.
We have $\omega_X=\pi^*(\omega_{\PP^2})\otimes\O_X((p-1)\Delta)$ by \cite{Rudakov; Safarevic 1976}, Proposition 2 combined with
Proposition 3, which gives $\O_X(\Delta)=\pi^*(\shL)$. Consequently
$$
H^0(X,\O_X(\Delta)) = H^0(\PP^2,\shA\otimes\shL) = H^0(\PP^2,\shL)\oplus H^0(\PP^2,\O_{\PP^2}) = k^n\times k,
$$
with the integer $n=(d+1)(d+2)/2$, as desired.

It is not difficult to compute bracket and $p$-map for $\lieh = H^0(X,\Theta_{X/k})$.
The coordinate rings for the affine open sets $U_i=D_+(T_i)$ of $\PP^2$ are the homogeneous localizations  $R_i=k[T_0,T_1,T_2]_{(T_i)}$,
and the preimages $\pi^{-1}(U_i)$ are the spectra of  $A_i=R_i[t_i]/(t^p_i-s_i)$. Here  $s_i$ denotes the dehomogenization of \eqref{defining relation}
with respect to $T_i$. We have $\Theta_{A_i/R_i}=A_i\partial/\partial t_i$, and 
our multiplicative vector field restricts becomes $D= t_i\partial/\partial t_i$.
For any $b_i,b'_i\in R_i$, one immediately calculates
$$
[b_i\partial/\partial t_i, t_i\partial/\partial t_i]=b_i\partial/\partial t_i,\quad
[b_i\partial/\partial t_i, b'_i\partial/\partial t_i]=0\quadand 
(b_i\partial/\partial t_i)^\pmap =0.
$$
Choosing a basis for $ H^0(\PP^2,\shL)$, we infer  that the vector space 
 decomposition $\lieh=H^0(\PP^2,\shL)\oplus H^0(\PP^2,\O_{\PP^2})$    becomes a  semi-direct product structure $\lieh=k^n\rtimes\liegl_1(k)$
for the restricted Lie algebra.
 
It remains to  check that  $\Sing(X/k)$ is finite, and that all singularities are rational double points. For this we may assume that $k$ is algebraically closed.
In light of the symmetry in \eqref{defining relation}, it suffices to verify this on the preimage $V=\pi^{-1}(U)$ of  the open set  $U=D_+(T_0)$.
Setting $x=T_1/T_0$ and $y=T_2/T_0$, we see that $V$ has coordinate ring $A=k[x,y,t]/(f)$ with
$$
f=t^p- xy-x^{pd-2}y-xy^{pd-2}.  
$$
The singular locus comprises the common zeros of $f$ and the partial derivatives  
$\partial f/\partial x=-y+2x^{pd-3}y-y^{pd-2}=0$ and $\partial f/\partial y= -x-x^{pd-2}+2xy^{pd-3}=0$.
Clearly there are only finitely many singularities with $x=0$ or $y=0$.
For the remaining part of  $\Sing(X)$, it suffices to examine the system of polynomial equations
\begin{equation}
\label{simplified derivatives}
-1+2x^{pd-3}-y^{pd-3}=0\quadand  -1-x^{pd-3}+2y^{pd-3}=0.
\end{equation}
This  is a system of linear equations in the powers $x^{pd-3}$ and $y^{pd-3}$, and $x^{pd-3}=y^{pd-3}=1$ is one solution.
Using $p\neq 3$ we see that there are no other solutions.
 
It remains to verify that all singularities are rational double points.   It would be tedious and cumbersome to do this explicitly.
We resort to a trick  of independent interest, where we actually show that there are only rational double points of $A$-type:
By Proposition~\ref{double points 1} and Lemma~\ref{double points 2}, it suffices to verify that no singular point is a zero for the polynomial 
$ (\frac{\partial^2 f}{\partial x \partial y})^2 -\frac{\partial^2 f}{\partial x^2}\cdot\frac{\partial^2 f}{\partial y^2}$,
which gives the additional equation
\begin{equation}
\label{double partials}
 1+4x^{2pd-6}+4y^{2pd-6}-4x^{pd-3}-4y^{pd-3}-28(xy)^{pd-3}=0.
\end{equation}
Substituting   $x^{pd-3}=y^{pd-3}=1$, the left-hand side becomes $1+4+4-4-4-28=-27$, which is indeed  non-zero because $p\neq 3$.
\qed

\medskip
Note that  $\pi:X\ra\PP^2$ is a universal homeomorphism, so the geometric fundamental group
of $X$ is trivial, and  $b_1=0$ and $b_2=1$.  
Let $x_1,\ldots,x_r$ be the  geometric singularities   on $X$. We saw  above 
that they are   rational double points  of certain $A_{n_i}$. One referee pointed out that
they all have $n_i=p-1$, which can be seen by considering the action of the Frobenius 
on the  local class group, which is multiplication by $p$ on a cyclic of order $n_i+1$.
Since the Frobenius factors over the projective plane, one infers that $n_i+1\mid p$,
hence $n_i+1=p$.
As discussed in Section \ref{Canonically polarized},
the Chern number $c_2$ is the alternating sum of the Betti numbers on the minimal resolution of $X$, which yields
the formula  $c_2-3=r(p-1)$. 

One referee also kindly pointed out that the arguments in the proof for Proposition \ref{example 1} hold true 
for \emph{general} polynomials $s\in \Gamma(\PP^2,\shL^{\otimes p})$ of degree $pd$,
provided $pd-d-2>0$ and $d\geq 2$, by using the result of Liedtke \cite{Liedtke 2013a}, Theorem 3.4, which 
ensures that all occurring  singularities must be  rational double points of type $A_{p-1}$.

Let us remark that  the surface $X\subset\PP^3$ defined by the homogeneous polynomial
$$
s=T_0T_1T_2^{2p-1}-T_0T_1T_3^{2p-1} + T_0^{2p} T_2 + T_1^{2p}T_3 +T_2^{2p+1} +T_3^{2p+1} 
$$
is a canonically polarized surface with $c_1^2=(2p-3)^2(2p+1)$ and $c_2=8p^3-4p^2+2p+3$, such that $\lieh=H^0(X,\Theta_{X/k})$
is isomorphic to $\liegl_1(k)$. We leave the details to the reader.

Next we   construct examples of smooth surfaces of general type $X$ where the restricted Lie algebra $\lieh=H^0(X,\Theta_{X/k})$
is isomorphic to $k^m$.
The possibility of the following construction was   suggested  by one of the referees:
Let $C$ be a smooth curve   with $h^0(\O_C)=1$, 
together with an   isomorphism $\varphi:\shL^{\otimes pl}\ra\Omega^1_{C/k}$
that is locally exact, for some invertible sheaf $\shL$ of degree $d\geq 1$,
and some integer $l\geq 1$ prime to the characteristic $p>0$.
This datum is called a \emph{generalized Tango curve} of index $l$. 

We are mainly interested
in the case $l\equiv -1$ modulo $p$, and then write $l=pn-1$.
Lang \cite{Lang 1983} used this situation  to construct smooth surfaces $X$ endowed with a fibration $f:X\ra C$ with $\O_C=f_*(\O_X)$,
where all geometric fibers are singular rational curves with a unique cuspidal singularity.
One also says that $(C,\shL,\varphi)$ is a \emph{generalized Tango curve} of type $(p,n,d)$,
and $X$ is the resulting \emph{generalized Raynaud surface} of type $(p,n,d)$.

\begin{proposition}\label{Raynaud surfaces}
\mylabel{example 3}
In the above setting, suppose that $p\geq 3$ and $n\geq 2$.  Then  $X$ is a minimal surface of general type with  
$$
\lieh= k^m \quadand \mathrm{Aut}_{X/k}[F] = \alpha^{\oplus m}_p,
$$
with $m=h^0(\shL)$. Moreover, the Chern invariants are given by the formulas  $c_1^2=d(p^4n^2 + 4p + 2np -n^2p^2-4np^2-2np^3)$
and $c_2= 2pd(1-np)$.
\end{proposition}

\proof
We may assume that $k$ is algebraically closed. Since we assume that our Tango curve has index $l\equiv -1$ modulo $p$, and the characteristic is $p\geq 3$,
we have $m=h^0(\shL)$, according to \cite{Takeda 1992}, Theorem 2.1.
Lang computed the Chern invariants, and observed that $X$ is minimal and of general type
(\cite{Lang 1983}, Theorem 2 and beginning of Section 2).

Both restricted  Lie algebras $k^m\rtimes \liegl_1(k)$ and $\liesl_2(k)$ contain non-zero multiplicative elements.
In light of Theorem \ref{structure result}, our task is 
to verify that non-zero multiplicative vector fields do not exist on $X$.
Seeking a contradiction, we suppose that   $\delta\in H^0(X,\Theta_{X/k})$ is such a vector field.
The  saturation  $\O_Y(\Delta)$ for the injection $\delta:\O_Y\ra\Theta_{Y/k}$ defines an effective Cartier divisor $\Delta\subset X$.
It follows from  \cite{Rudakov; Safarevic 1976}, Theorem 2 that every connected component is smooth. Hence 
the   irreducible components are pairwise disjoint, and  horizontal for the fibration $f:X\ra C$,
because all closed fibers are singular.

To reach a contradiction we examine various curves on $X$ and their intersection numbers.
Write $F\subset X$ for a closed fiber, $D\subset X$ for reduced support of $\Omega^1_{X/C}$, and $S\subset X$ for
the canonical section  constructed in \cite{Lang 1983}, Section 2. Note that $D$ is also called the curve of cusps.
According to loc.\ cit.\ one has 
$$
S^2=d,\quad F^2=0,\quad (F\cdot S)=1, \quadand D=-pdF+pS.
$$
Consequently $D^2=-p^2d$. For the curve $G=dF+nD$ we get $G^2=dnp(2-np)<0$.
According to loc.\ cit., Theorem 1 there is an exact sequence
$$
0 \lra \O_X(G)\lra \Theta_{X/k} \lra\omega_X^{\otimes-1}(-G) \lra 0,
$$
giving an identification $H^0(X,\O_X(G))=H^0(X,\Theta_{X/k})$. Our global vector field $\delta$
factors over the inclusion $\O_X(G)$, which gives an equality $\O_X(\Delta)=\O_X(G)$ as subsheaves of $\Theta_{X/k}$.
In particular, the curves $\Delta$ and $G$ are linearly equivalent.
Decompose the smooth curve $\Delta=\Delta_1+\ldots+\Delta_r$ into irreducible components. 
We already observed that each $\Delta_i$ is  horizontal.
From  $\Delta\cdot (dF+nD)=G^2<0$ we infer $D\subset\Delta$. Now consider   $\Delta'=\Delta-D$,
which contains neither $D$ nor $F$, and   is linearly equivalent to $G'=dF+(n-1)D$.
Then
$$
pd(n-1)(2-(n-1)p)=(G')^2=(dF+(n-1)D)\cdot \Delta'\geq 0.
$$
By our assumptions we have $p\geq 3$ and $n\geq 2$, hence the left hand side is strictly negative, contradiction.
\qed

\medskip
There are indeed generalized Tango curves $C$ of type $(p,n,d)$ with non-zero $m=h^0(\shL)$,
for instance the curve $C$ with affine equation $y^{lp}-y=x^{lp-1}$, where we set $l=pn-1$,
according to \cite{Takeda 1992}, Example 1.2.

It remains to  verify some  technical results  used throughout this section. 
The following results are   well-known over the field of complex numbers 
(compare for example  \cite{Kollar; Mori 1998}, Section 4.2 and \cite{Arnold 1985}, Part II). 
The  arguments  apparently work in all characteristics except $p=2$. 
For the convenience of the reader we give   self-contained and characteristic-free proofs.
 
Let  $A$ be a complete local $k$-algebra  
that is regular of dimension three, with maximal ideal $\maxid_A$ and residue field $k=A/\maxid_A$. 
Note that for each  choice of regular   system of parameters
$x,y,z\in A$ one obtains an identification $A=k[[x,y,z]]$.
 
\begin{lemma}
\mylabel{double points 1}
Let $f\in A$ be an irreducible element such that  $f\equiv x_0y_0+\lambda z_0^2$ modulo $\maxid_A^3$ for some regular system of parameters
$x_0,y_0,z_0$ and some  $\lambda\in k$. Then 
there exists another regular system of parameters $x,y,z$   such that $(f)=(xy+z^n)$, for some $n\geq 2$.
Hence $B=A/(f)$ is a rational double point   of type $A_{n-1}$. 
\end{lemma}

\proof
We construct by induction on $n\geq 0$ certain  regular system of parameters  $x_n,y_n,z_n\in A$ such  that 
$x_n\equiv x_{n-1}$ and  $y_n\equiv y_{n-1}$ modulo $\maxid_A^n$, and $z_n=z_0$, and 
\begin{equation}
\label{initial description}
f=x_ny_n +x_n\phi_n+y_n\psi_n+h_n
\end{equation}
for some $\phi_n, \psi_n \in \mathfrak{m}_A^{n+2}$ and $h_n\in z_n^2k[[z_n]]$.
For $n=0$ we   take $x_0,y_0, z_0$ as in our assumptions.
If we already have defined $x_n,y_n,z_n\in A$, we set
\begin{equation}
\label{inductive change}
x_{n+1} = x_n + \psi_n,\quad y_{n+1}=y_n+\phi_n\quadand z_{n+1}=z_n.
\end{equation}
Clearly   $x_{n+1}\equiv  x_n$ and $y_{n+1}\equiv  y_n$ modulo $\maxid_A^{n+1}$, and $z_{n+1}=z_0$.
In particular the above is a regular system of parameters.
Since $\phi_n\psi_n\in\maxid_A^{2n+4}$, we may write 
$$
-\phi_n\psi_n = x_{n+1}\phi_{n+1} + y_{n+1}\psi_{n+1} + \sigma_n 
$$
with $\phi_{n+1},\psi_{n+1}\in \maxid_A^{2n+3}$ and $\sigma_n\in z_{n+1}^{2n+4}k[[z_{n+1}]]$. 
Combining \eqref{initial description} and \eqref{inductive change}, 
we get
$$
f=x_{n+1}y_{n+1} +x_{n+1}\phi_{n+1} + y_{n+1}\psi_{n+1} + h_{n+1}
$$
where $h_{n+1}=h_n +\sigma_n$ belongs to $z_{n+1}^2k[[z_{n+1}]]$. This completes our inductive definition.
Note that $h_{n+1}\equiv h_n$ modulo $\maxid_A^{2n+4}$.

By construction, the $x_n,y_n,z_n$ are convergent sequences in $A$ with respect to the $\maxid_A$-adic topology.
The limits $x,y,z\in A$ give the desired regular system of parameters: Since the $\phi_n,\psi_n$ converge to zero,
we have $f=xy+h(z)$, where $h$ is the limit of the $h_n\in k[[z]]$.
We must have $h\neq 0$, because $f$ is irreducible. Hence $h=uz^n$ with  $u\in k[[z]]^\times$  
and $n\geq 2$. Replacing $x$ by $u^{-1}x$, we finally get $(f)=(xy-z^n)$. Summing up,  $B=A/(f)$
is  a rational double point of type $A_{n-1}$.
\qed

\medskip
The condition in the proposition  can be checked with partial derivatives, at least   
if $k$ is algebraically closed. This makes the criterion applicable for computations:

\begin{proposition}
\mylabel{double points 2}
Let $f\in\maxid_A^2$.  Suppose that $k$ is algebraically closed, and that 
\begin{equation}
\label{second partials}
\left(\frac{\partial^2 f}{\partial u_1 \partial u_2}\right)^2-
\left(\frac{\partial^2 f}{\partial u_1^2}\right)\cdot\left(\frac{\partial^2 f}{\partial u_2^2}\right)
\not\in \mathfrak{m}_A 
\end{equation}
for some system of parameters $u_1,u_2,u_3\in A$.
Then there exists another system of parameters $x,y,z$ such that $f\equiv xy+\lambda z^2$ modulo $\mathfrak{m}_A^3$, for some $\lambda \in k$.
\end{proposition}

\proof
Write  $f=q+g$, where  $q=q(u_1,u_2,u_3)$ is a homogeneous polynomial of degree two 
and $ g\in \maxid_A^3$.  Write $q=q_1 +u_3l $, where $l=l(u_1,u_2,u_3)$ is homogeneous of degree one
and $q_1=q_1(u_1,u_2)$. 
If $q_1$ is a square, a straightforward computation with   partial derivatives produces
a contradiction to \eqref{second partials}.
Since $k$ is algebraically closed
we have a factorization $q_1=L_1\cdot L_2$
where   $L_1=L_1(u_1,u_2)$ and $L_2= L_2(u_1,u_2)$ are independent 
homogeneous polynomials of degree  one.
Then $w_1=L_1$, $w_2=L_2$ and $w_3$ form another regular system of parameters of $A$, and we have 
$$
q=w_1w_2+w_3l =w_1w_2+aw_3w_1+bw_3w_2+cw_3^2,
$$
with $a,b,c \in k$. 
We finally set  $x=w_1+bw_3$, $y=w_2+aw_3$ and $z=w_3$. This is a further
regular  system of parameters, with $q=xy+\lambda z^2$, where $\lambda=c-ab$.  Therefore, $f\equiv xy+\lambda z^2$ modulo $\mathfrak{m}_A^3$, as claimed.
\qed

\medskip
We also used a general fact on coherent sheaves: Let $X$ be a noetherian   scheme that is integral and normal,
$\shE$ be a coherent sheaf of rank two, 
$s:\O_X\ra\shE^\vee$ a non-zero global section. The double dual $\O_X(-\Delta)$ for the image
of the dual map $s^\vee:\shE^{\vee\vee}\ra \O_X$ defines an effective Weil divisor $\Delta\subset X$.
By \cite{Hartshorne 1994}, Corollary 1.8 the duals of coherent sheaves on $X$ are reflexive. Dualizing 
$\shE^{\vee\vee}\ra\O_X(-\Delta)\subset  \O_X$ we see that the homomorphism $s$ factors over an inclusion $\O_X(\Delta)\subset \shE^\vee$.
The latter is called the  \emph{saturation} of the section $s\in\Gamma(X,\shE^\vee)$.

\begin{lemma}
\mylabel{four term sequence}
In the above setting  there is a four-term exact sequence
$$
0\lra \O_X(\Delta)\lra \shE^\vee\lra \shL(-\Delta) \lra \shN\lra 0,
$$
where $\shL=\uHom(\Lambda^2(\shE),\O_X)$, and $\shN$ is a coherent sheaf whose support has codimension at least two.
\end{lemma}

\proof
According to \cite{Hartshorne 1994}, Theorem 1.12 it suffices to construct a short exact sequence $0\ra\O_X(\Delta)\ra\shE^\vee\ra\shL(-\Delta)\ra 0$
on the complement of some closed set $Z\subset X$ of codimension at least two.
Let $\shE_0$ be the quotient of $\shE$ by its torsion subsheaf. The surjection $\shE\ra\shE_0$ induces
an equality $\shE_0^\vee=\shE^\vee$, so we may assume that $\shE$ is torsion free. It is then locally free in codimension one,
so it suffices to treat the case that $\shE$ is locally free. By construction,
the cokernel $\shF$ for $\O_X(\Delta)\subset\shE^\vee$ is torsion-free of rank one, so we may assume that  it is invertible.
Taking determinants shows $\shF\simeq \shL(-\Delta)$.
\qed


\end{document}